# ORTHOMODULARITY IN INFINITE DIMENSIONS; A THEOREM OF M. SOLÈR

SAMUEL S. HOLLAND, JR.

ABSTRACT. Maria Pia Solèr has recently proved that an orthomodular form that has an infinite orthonormal sequence is real, complex, or quaternionic Hilbert space. This paper provides an exposition of her result, and describes its consequences for Baer ∗-rings, infinite-dimensional projective geometries, orthomodular lattices, and Mackey's quantum logic.

## 1. The theorem of Maria Pia Solèr

This theorem deals with infinite-dimensional Hermitian forms over ∗-fields. A ∗-field $\mathcal{K}$ is a (not necessarily commutative) field with involution, an involution being a map $\rho \mapsto \rho^*$ of $\mathcal{K}$ onto itself that satisfies $(\rho + \sigma)^* = \rho^* + \sigma^*$, $(\rho\sigma)^* = \sigma^*\rho^*$, and $\rho^{**} = \rho$ for all $\rho$, $\sigma$ in $\mathcal{K}$. The ∗-fields that appear in the conclusion of Solèr's Theorem are the three classical examples: the real number $\mathbb{R}$ with the identity involution; the complex numbers $\mathbb{C}$ with complex conjugation as involution, $(\zeta + i\eta)^* = \zeta - i\eta$; and the real quaternions $\mathbb{H}$ consisting of all $\zeta_0 + \zeta_1 i + \zeta_2 j + \zeta_3 k$, $\zeta_i \in \mathbb{R}$, multiplication determined by $i^2 = j^2 = -1$ and $ij = -ji = k$, and with quaternionic conjugation as involution, $(\zeta_0 + \zeta_1 i + \zeta_2 j + \zeta_3 k)^* = \zeta_0 - \zeta_1 i - \zeta_2 j - \zeta_3 k$.

Contrary to what our early graduate education would lead us to believe, ∗-fields exist in incredible abundance and variety; the three classical number fields of analysis, described above, are extraordinarily special within the class of general ∗-fields. All the more startling then that the theorem of Solèr, which, on the basis of two seemingly innocent assumptions, namely, that the form is orthomodular and has an infinite orthonormal sequence (terms we shall come to shortly), concludes that the underlying ∗-field is in fact either $\mathbb{R}$, $\mathbb{C}$, or $\mathbb{H}$.

We set the stage for a statement of her theorem.

**Definition 1.1.** Consider a left vector space $E$ of any dimension, finite or infinite, over a ∗-field $\mathcal{K}$. A **Hermitian form** $<\cdot,\cdot>$ on $E$ is a mapping $E \times E \longrightarrow \mathcal{K}$ that associates to every pair of vectors $x$, $y$ in $E$ a scalar $<x,y>$ in $\mathcal{K}$ in accordance with the following rules:

(1) $\begin{cases} <\rho x + \sigma y, z> = \rho <x,z> + \sigma <y,z> \\ <x, \rho y + \sigma z> = <x,y>\rho^* + <x,z>\sigma^* \end{cases}$ for all $x$, $y$, $z$ in $E$ and all $\rho$, $\sigma$ in $\mathcal{K}$.

(2) If either $<a,x> = 0$ $\forall x \in E$, or $<x,a> = 0$ $\forall x \in E$, then $a = 0$.

(3) $<x,y>^* = <y,x>$ for all $x$, $y$ in $E$.

A function that satisfies just (1) and (2) is called variously "conjugate-bilinear" or "sesquilinear". Rule (2), expressing nonsingularity, is not normally part of the definition of a form, but I have included it here for convenience. Rule (3) makes the form Hermitian; the term "symmetric" is sometimes used when $*=$identity.

Suppose given a $\mathcal{K}$-space $E$ with Hermitian form $<\cdot,\cdot>$. Two vectors $x$, $y$ in $E$ are **orthogonal** when $<x,y>=0$. A sequence $\{e_i : i = 1, 2, ...\}$ of nonzero vectors in $E$ is called orthogonal when $<e_i, e_j>=0$ for $i \neq j$ and is called **orthonormal** when also $<e_i, e_i>=1$ for all $i = 1, 2, ...$.

Given a nonempty subset $S$ of $E$, the symbol $S^\perp$ (sometimes read "$S$-perp") stands for the set of those elements of $E$ orthogonal to every element of $S$: $S^\perp = \{x \in E : <x,s>=0 \ \forall s \in S\}$. $S^\perp$ is always a subspace of $E$ irrespective of the nature of the nonempty subset $S$.

A subspace $M$ of $E$ is called **closed** when $M = M^{\perp\perp}$. Every finite-dimensional subspace of $E$ is closed. Let us prove this for a one-dimensional subspace $M = \mathcal{K}x$, $0 \neq x \in E$. As our form is nonsingular (axiom (2)), there is an $a \in E$ such that $<a,x> \neq 0$. Then for any $u \in E$, the vector $z = u - <u,x><a,x>^{-1} a$ is orthogonal to $x$:

$$<z, x> = <u - <u, x><a, x>^{-1} a, x>$$
$$= <u, x> - <u, x> = 0.$$

Hence $z \in (\mathcal{K}x)^\perp$ for every $u \in E$. So, given $y \in (\mathcal{K}x)^{\perp\perp}$, we must have $<z, y> = 0$ for every $u$; $<u - <u, x><a, x>^{-1} a, y> = 0$ for every $u \in E$. This equation can be rewritten $<u, y - \rho^* x> = 0$ for every $u \in E$, where $\rho = <a, x>^{-1} <a, y>$. Hence $y - \rho^* x = 0$, so $y = \rho^* x \in \mathcal{K}x$. Thus $(\mathcal{K}x)^{\perp\perp} \subseteq \mathcal{K}x$ and, as the inclusion $M \subseteq M^{\perp\perp}$ always holds, we have $(\mathcal{K}x)^{\perp\perp} = \mathcal{K}x$. Using the same method, one can show that for any subspace $M$ of $E$ and any $x \in E$, $(M + \mathcal{K}x)^{\perp\perp} = M^{\perp\perp} + \mathcal{K}x$. In particular if $M$ is closed and $N$ is finite-dimensional, the $M + N$ is closed.

The Hermitian form is called **orthomodular** when $M + M^\perp = E$ for every closed subspace $M$. We repeat

**Definition 1.2.** The Hermitian space $\{E, \mathcal{K}, <\cdot,\cdot>\}$ is **orthomodular** if, and only if,
$$\emptyset \neq M \subseteq E \ \& \ M = M^{\perp\perp} \implies M + M^\perp = E.$$

A Hermitian form is called **anisotropic** when $x \neq 0 \implies <x,x> \neq 0$. An orthomodular form is necessarily anisotropic because were there a nonzero $x$ with $<x,x>=0$, then $\mathcal{K}x \subseteq (\mathcal{K}x)^\perp$, so $(\mathcal{K}x) + (\mathcal{K}x)^\perp = (\mathcal{K}x)^\perp \neq E$ (if $(\mathcal{K}x)^\perp = E$, then $\mathcal{K}x = (\mathcal{K}x)^{\perp\perp} = E^\perp = 0$, a contradiction). In the finite-dimensional case, the converse also holds: anisotropic $\implies$ orthomodular. Because anisotropic implies that $M \cap M^\perp = 0$ for any subspace $M$, so $(M + M^\perp)^\perp = M^\perp \cap M = 0$. Thus $(M + M^\perp)^{\perp\perp} = 0^\perp = E$. But $M + M^\perp$ is closed, so $M + M^\perp = E$.

Thus in the finite-dimensional case "anisotropic" and "orthomodular" mean the same thing. Passing to the infinite-dimensional case, we seem to have a choice: generalize as anisotropic, or generalize as orthomodular. Requiring an infinite-dimensional form to be anisotropic turns out to be a feeble restriction. All the more dramatic then the **1.3   Theorem of Maria Pia Solèr** [S]. *Let $\mathcal{K}$ be a $*$-field, $E$ a left vector space over $\mathcal{K}$, and $<\cdot,\cdot>$ an orthomodular form on $E$ that has an infinite orthonormal sequence. Then $\mathcal{K} = \mathbb{R}, \mathbb{C},$ or $\mathbb{H}$, and $\{E, \mathcal{K}, <\cdot,\cdot>\}$ is the corresponding Hilbert space.*

Clearly both assumptions, orthomodularity and the existence of an infinite orthonormal sequence, are **necessary** if one wishes to characterize infinite-dimensional Hilbert space, because both are true there. Orthomodularity is just the well-known projection theorem, and the existence of an infinite orthonormal sequence comes out of the basic geometry plus the fact that positive real numbers have square roots. The fact that these two conditions are also **sufficient** has to be one of the more surprising results of infinite-dimensional algebra, even mathematics generally, because there is absolutely nothing in either condition to suggest the presence of the three number fields of analysis. Also, no amount of computer simulation would suggest this result.

One may relax the assumption that $\{E, \mathcal{K}, <\cdot,\cdot>\}$ has an infinite orthonormal sequence to the assumption that $E$ contains an orthogonal sequence $\{e_i : i \in \mathbb{N}\}$ of nonzero vectors such that $<e_i, e_i> = <e_j, e_j>$ for all $i, j$. Because if that is the case, then setting $\lambda = <e_i, e_i>$, we define a new involution $\#$ on $\mathcal{K}$ by $\rho^\# = \lambda \rho^* \lambda^{-1}$ and a new form $[\cdot,\cdot]$ on $E$ by $[\cdot,\cdot] = <\cdot,\cdot> \lambda^{-1}$. A direct calculation shows that the new form $[\cdot,\cdot]$ is Hermitian with respect to the new involution $\#$ and is also orthomodular because it induces the same $\perp$ map. The sequence $\{e_i\}$ is now orthonormal in $\{E, \mathcal{K}(\#), [\cdot,\cdot]\}$, so by Solèr's Theorem $\mathcal{K}(\#)$ is $\mathbb{R}, \mathbb{C},$ or $\mathbb{H}$. But $\lambda = \lambda^\#$, so $\lambda$ is a nonzero real number. Thus $\# = *$, so $\mathcal{K}(*)$, our original $*$-field, is $\mathbb{R}, \mathbb{C},$ or $\mathbb{H}$, and either $<\cdot,\cdot>$ or $-<\cdot,\cdot>$ is a Hilbert space inner product.

If we drop the assumption that the quadratic form assumes a common value on an infinite orthogonal sequence and require only that $\{E, \mathcal{K}, <\cdot,\cdot>\}$ is orthomodular and infinite-dimensional, what then can be said? Very little is known.

Solèr's Theorem has a history. Kaplansky initiated the study of infinite-dimensional forms with his 1950 paper [Kap1]. In this paper he proved that an infinite-dimensional orthomodular form cannot have countable Hamel dimension [Kap1, p. 4], which was the first use of the combination orthomodularity and infinite dimension.

In 1964, Piron published a theorem asserting that an orthomodular inner product space over $\mathbb{C}$ was metrically complete [Pir, Theorem 22]. While Piron's proof contained an error, he is to be credited with surmising the result which was subsequently proved by Amemiya and Araki [A-A]. In the acknowledgements to their paper, Araki thanks Marshall Stone as well as Piron for bringing the problem to his attention. The Amemiya-Araki method of proof, which works just as well in the real and quaternionic cases, is nicely detailed in Maeda's book [M-M, Theorem 34.9]. It has found frequent application since its publication.

The Amemiya-Araki-Piron Theorem was extended by various people over the period 1970–1977 to cover the case $\mathcal{K}$ is a $*$-closed subfield of $\mathbb{H}$. This essentially took care of the case where $\mathcal{K}$ had an archimedean Baer ordering. At this point the strength of the orthomodular axiom began to be more apparent, which gave rise to speculation that orthomodularity by itself would characterize Hilbert space. In 1980, Hans A. Keller put a decisive end to this speculation [Ke]. Keller constructed a "nonclassical Hilbert space"; i.e. an infinite-dimensional orthomodular form over a $*$-field $\mathcal{K}$ different from $\mathbb{R}$, $\mathbb{C}$, or $\mathbb{H}$. Of course, his form did not have an infinite orthonormal sequence. Keller's construction sparked much interest and was much generalized. Valuation theory emerged as a key tool in all these generalizations.

Keller's result showed that orthomodularity by itself did not suffice, but, as remarked above, adding in an archimedean ordering together with an infinite orthonormal sequence tipped the scales.

What about a nonarchimedean ordering? How to get rid of the infinitesimal elements? The answer came in a brilliant 1977 paper by W. John Wilbur [W]. Wilbur proved that an orthomodular space $\{E, \mathcal{K}, <\cdot, \cdot>\}$ is Hilbert space provided that the $*$-field $\mathcal{K}$ is either commutative or 4-dimensional over its center and has the property that for every $\rho = \rho^* \in \mathcal{K}$ there is an $\alpha \in \mathcal{K}$ such that $\pm\rho = \alpha\alpha^*$ [W, Theorem 5.8]. (Under this assumption the $*$-field $\mathcal{K}$ is uniquely ordered, and the form $<\cdot, \cdot>$ must represent 1 on every one-dimensional subspace. In [H3] Wilbur's theorem was generalized to Baer ordered $*$-fields that have this Henselian property: given any infinitesimal $\varepsilon = \varepsilon^*$ in $\mathcal{K}$, there is an $\alpha \in \mathcal{K}$ such that $1 + \varepsilon = \alpha\alpha^*$.) While Wilbur did put some stringent conditions on the $*$-field $\mathcal{K}$, he introduced no assumptions, tacit or explicit, about $\mathcal{K}$ that would imply that its ordering was archimedean. To the contrary, he devised very clever and original arguments to eliminate the infinitesimal elements, and this was a major breakthrough.

Wilbur's impact on the final resolution of this problem goes further. His ideas, as subsequently extended and generalized by Gross and Keller, constitute key ingredients in Solèr's proof. In his paper [W] Wilbur proved that, given two mutually orthogonal orthogonal sequences, $\{e_i : i = 1, 2, ...\}$ and $\{f_i : i = 1, 2, ...\}$, in an orthomodular space $\{E, \mathcal{K}, <\cdot, \cdot>\}$, and given $x \in \{e_i : i = 1, 2, ...\}^{\perp\perp}$, then there exists $y \in \{f_i : i = 1, 2, ...\}^{\perp\perp}$ such that $<y, f_i> = <x, e_i>$, $i = 1, 2, ...$, and $<y, y> = <x, x>$. This result was generalized in several directions by Gross [G2, §I.1]. Gross's work, in turn, forms an essential component in Solèr's proof [S, Lemma 3 and Corollaries 1–7]. Another one of Wilbur's results [W, Lemma 5.3] was substantially extended by Keller to show that $\mathbb{R} \subseteq \mathcal{K}$ and (roughly) $\ell_2(\mathbb{R}) \subseteq E$ [G2, Theorem I.2.1]. Keller's result was further generalized by Solèr and appears as Theorem 1 in her paper [S]. Wilbur's paper [W] was a watershed event, and it stands out as the signal contribution prior to the decisive result of Solèr.

Finally, one name needs to be placed in special prominence in connection with this theorem—the name of Herbert Gross. Gross greatly influenced the development of the theory of infinite-dimensional forms, both by his many papers and his book [G1], and especially by his continuous stream of excellent students. His productivity was all the more admirable, given his never-ending battle against a mysterious refractory illness. One of his last papers [G2] is a beautiful study of the orthomodular problem which was finally solved by Solèr. Both Hans A. Keller and Maria Pia Solèr were students of Herbert Gross. Upon his sudden death on 29 October 1989, Solèr, who had just begun work on Gross's orthomodular problem, began to look for someone else to guide her thesis. She came in contact with Professor A. Prestel (at Konstanz), a friend of Gross, who encouraged her to continue work on the problem and offered to guide her thesis unofficially. Work began under that arrangement, with Professor H. Storrer at the University of Zurich as her offical advisor. Prestel suggested some ideas, but very soon Solèr began working independently. While she had some precious well-appreciated help from a very interested Hans Keller, who had just returned to Zurich, she otherwise worked alone to prove the result we are expounding here.

Other proofs of Solèr's Theorem have now appeared [KKS, Prs].

## 2. Baer ∗-rings

Kaplansky's Theory of Baer ∗-rings provides an algebraic setting for the main features of the von Neumann-Murray theory of operator algebras. The motivations and history of the subject are nicely set forth in the preface to Kaplansky's 1968 book [Kap2]. (The work itself dates from the early 1950's.) A **Baer ∗-ring** A is a ring with identity that carries an involution $x \mapsto x^*$ and has the property that the right annihilator of any nonempty subset of $A$ is a principal right ideal generated by a projection. A **projection** $e \in A$ is a symmetric idempotent: $e^2 = e^* = e$. In more detail: given $\emptyset \neq S \subseteq A$, set $RtAnn(S) = \{x \in A : sx = 0 \; \forall s \in S\}$. Then the condition defining a Baer ∗-ring is this: $\emptyset \neq S \subseteq A \Longrightarrow RtAnn(S) = eA$ for some projection $e \in A$. If we apply the involution, we get the same statement with "right" replaced by "left", so it does not matter which we use. A **Baer ∗-factor** is a Baer ∗-ring whose only central projections are $0$ and $1$. Factors display all the essential characteristics of the theory, and we shall concentrate on them.

Kaplansky proved straight off that, just as in the theory of von Neumann algebras, factors come in five kinds. These are described by combining the familiar terms "types I, II, and III" with the words "finite" and "infinite". As our only concern here is with infinite type I Baer ∗-factors, we shall describe in detail only these. The description is in terms of the behavior of the set of projections within the factor. These projections form a partially ordered set under the ordering $e \leq f \iff ef = e$; there is also an orthogonality relation: $e$ and $f$ are orthogonal when $ef = 0$ (equivalently $e \leq 1 - f$). A nonzero projection $e$ in a factor $A$ is **minimal** if $f \leq e \Longrightarrow f = 0$ or $f = e$, and $A$ is **type I** if it has a minimal projection. I call $A$ **infinite** if it contains an infinite family of nonzero mutually orthogonal projections, and any family of nonzero mutually orthogonal projections is at most countable. The latter condition is usually called "separable", but I lump it in here for convenience.

Two projections $e$, $f$ in $A$ are said to be **equivalent relative to** $A$ if there is a $w \in A$ such that $w^*w = e$ and $ww^* = f$. This key notion of relative equivalence (often called just equivalence) is taken over word for word and letter for letter from the theory of von Neumann algebras. But taking over to our purely algebraic environment the body of theory based on this notion of equivalence is another matter. The simple axiom for a Baer ∗-ring will not support a theory of equivalence corresponding to that in von Neumann algebras. To make possible a purely algebraic theory of equivalence in Baer ∗-rings, Kaplansky introduced the **Existence of Projections** (EP) axiom: For every nonzero $x \in A$ there is a $y = y^*$ in $A$ such that $y$ commutes with everything that commutes with $xx^*$ and $xx^*y^2$ is a nonzero projection [Kap2, §13].

An infinite type I Baer ∗-factor that satisfies the EP axiom possesses a theory of relative equivalence parallel in every respect to that of its analytic prototype, the algebra of all bounded linear operators on infinite-dimensional separable complex

Hilbert space. One may consult Berberian's richly detailed book [Be] for a beautiful development of this purely algebraic theory.

The success of this algebraic version of the von Neumann-Murray theory brought up the question—are there actual examples of nonclassical $*$-rings to which it applies? Over the years no such examples appeared; there seemed to be no other examples of infinite type I Baer $*$-factors satisfying EP other than the classical examples. This fact prompted the enunciation, in 1973, of the following project [H2, p. 522]:

> *Find all infinite type* I *Baer $*$-factors that satisfy the EP axiom.*

Now, twenty years later, we have found them all:

**Theorem 2.1.** *An infinite type* I *Baer $*$-factor that satisfies the EP axiom is the algebra of all bounded linear operators on an infinite-dimensional separable real, complex, or quaternionic Hilbert space.*

This theorem follows from that of Solèr. Let $A$ be the Baer $*$-factor in question. As $A$ is type I, it has a minimal projection $e_0$. Then $\mathcal{K} = e_0 A e_0$ is a $*$-field under the involution inherited from $A$. And, under the addition and multiplication of $A$, the set $E = e_0 A$ is a left vector space over $\mathcal{K}$. Define a Hermitian form on $E$ as follows: Given $x = e_0 a$, $y = e_0 b$ in $E$, set $<x, y> = (e_0 a)(e_0 b)^* = e_0 a b^* e_0 \in \mathcal{K}$. This form is anisotropic because in a Baer $*$-ring $xx^* = 0 \Longrightarrow x = 0$. We shall apply Solèr's Theorem to this Hermitian space $\{E, \mathcal{K}, <\cdot, \cdot>\}$, and to do that, we need to verify the two hypotheses: (1) orthomodular and (2) existence of an infinite orthonormal sequence.

Construct the ring $B(E, \mathcal{K})$ of all $\mathcal{K}$-linear everywhere-defined operators on $E$ that have adjoints with respect to our form. Denote the adjoint also by $*$. Hence $B(E, \mathcal{K})$ consists of all those everywhere-defined $\mathcal{K}$-linear transformations $t$ of $E$ into itself for which there exists a (necessarily unique) $t^* \in B(E, \mathcal{K})$ such that $<(x)t, y> = <x, (y)t^*>$ for all $x, y$, in $E$ (writing the operators on the right). The operation of left multiplication by an element $e_0 z e_0 \in$ center$(\mathcal{K})$ is in $B(E, \mathcal{K})$, and the adjoint of this operation is left multiplication by $e_0 z^* e_0$.

Given $x \in A$, associate to it the operator $\hat{x}$ on $E$ given by $(e_0 a)\hat{x} = e_0 ax$. The operator $\hat{x}$ is $\mathcal{K}$-linear and has an adjoint: $(\hat{x})^* = (x^*)^\wedge$. The correspondence $x \mapsto \hat{x}$ is a $*$-isomorphism of the $*$-ring $A$ into the $*$-ring $B(E, \mathcal{K})$, so we may, and shall, identify $A$ with its image $\hat{A} \subseteq B(E, \mathcal{K})$. $A$ is also an algebra over center$(\mathcal{K})$, because left multiplication by $e_0 z e_0 \in$ center$(\mathcal{K})$ is the operator $(e_0 z e_0)^\wedge$.

Next observe that $A$ contains all the rank-one operators and all the projections in $B(E, \mathcal{K})$. Each rank-one operator $t \in B(E, \mathcal{K})$ has the form $(\cdot)t = <\cdot, y> z$ for some fixed $y, z \in E = e_0 A$. In our case we shall have $y = e_o \ell$ and $z = e_0 m$ for some $\ell, m \in A$, and a direct calculation shows that if $s = \ell^* e_0 m$, then $\hat{s} = t$ (because $(x)\hat{s} = e_0 a \ell_0^* e_0 m = (e_0 a)(e_0 \ell)^* e_o m = <x, y> z = (x)t$ for all $x = e_0 a \in A$). Hence $A$ contains all the rank-one operators in $B(E, \mathcal{K})$. Then the lemma from [H2] tells us that a subspace $M$ of $E$ has the form $M = (E)\hat{e}$ for some projection $e \in A \Longleftrightarrow M = M^{\perp\perp}$. Thus the closed subspaces of $E$ correspond one-to-one to the projections in $A$. But the same is true of the projections in $B(E, \mathcal{K})$; every projection $t \in B(E, \mathcal{K})$ is uniquely determined by its image, namely, the closed subspace $M$ given by $M = (E)t$. Thus every projection $t \in B(E, \mathcal{K})$ has the form $t = \hat{e}$ for some projection $e \in A$, so $A$ contains all the projections in $B(E, \mathcal{K})$.

The orthomodularity of our space $\{E, \mathcal{K}, <\cdot, \cdot>\}$ is now immediate, because given a closed subspace $M$ of $E$, there is a projection $e \in A$ such that $M = (E)\hat{e}$; hence $E = (E)\hat{e} + (E)(1-e)^{\wedge} = M + M^{\perp}$.

As for the required infinite orthonormal sequence, the EP axiom gives us more: given any $x \in E$, $x \neq 0$, there is $0 \neq \lambda = \lambda^* \in \mathcal{K}$ with $<\lambda x, \lambda x> = 1$. Which is to say, in any direction there is a vector of length 1. Because, given $x = e_0 x \in A$, the EP axiom tells us that there is a symmetric $y = y^* \in A$ that doubly commutes with $(e_0 x)(e_0 x)^* = e_0 x x^* e_0$ and makes $e_0 x x^* e_0 y^2 = p$ a nonzero projection. As $pe_0 = p$, and $e_0$ is minimal, we must have $p = e_0$. The projection $e_0$ commutes with $(e_0 x)(e_0 x)^*$, so $y$ commutes with $e_0$ and commutes with $(e_0 x)(e_0 x)^*$ as well. Hence we have

$$e_0 = (e_0 y e_0)(e_0 x x^* e_0)(e_0 y e_0).$$

Now $\lambda = e_0 y e_0$ is a symmetric element of our $*$-field $\mathcal{K}$, and $e_0 = e_0 1 e_0$ is the 1 of $\mathcal{K}$. Hence the above-displayed formula translates to $1 = \lambda <x, x> \lambda = <\lambda x, \lambda x>$ as desired. Our vector space $E$ has an infinite orthogonal sequence: take any sequence $\{e_i : i = 1, 2, ...\}$ of nonzero orthogonal projections in $A$, select one nonzero vector from each of the subspaces $(E)\hat{e}_i$, then, using the EP axiom, normalize each vector to length 1 by the above procedure.

Hence the space $\{E, \mathcal{K}, <\cdot, \cdot>\}$ is orthomodular and has an infinite orthonormal sequence. Solèr's theorem applies: $\mathcal{K} = \mathbb{R}, \mathbb{C}$ or $\mathbb{H}$, and $\{E, \mathcal{K}, <\cdot, \cdot>\}$ is the corresponding infinite-dimensional, separable, real, complex or quaternionic Hilbert space. The $*$-ring $B(E, \mathcal{K})$ is accordingly the full $*$-ring of all bounded linear operators on the Hilbert space $E$ and is an algebra over $\mathbb{R}, \mathbb{C}$, or $\mathbb{R}$ respectively. Our Baer $*$-ring $A$ is a $*$-subalgebra of $B(E, \mathcal{K})$ that contains all its projections; therefore, $A = B(E, \mathcal{K})$ [Dix, Proposition 7; F-T; H4]. Theorem 2.1 is proved.

What can be said about the coordinatization of Baer $*$-factors of other types? Nothing seems to be known.

## 3. Infinite-dimensional projective geometries

The fundamental theorem of projective geometry associates to each synthetic projective geometry $\mathbb{P}$ of dimension $\geq 3$ a field $\mathcal{K}$ which is uniquely determined up to isomorphism by $\mathbb{P}$. Solèr's Theorem furnishes a characterization of those infinite-dimensional geometries for which $\mathcal{K}$ is $\mathbb{R}, \mathbb{C}$, or $\mathbb{H}$.

Much of the development of this section has already been done very elegantly by S. Maeda [M–M, §§33, 34]. S. Maeda's exposition is couched in the language of lattice theory. Our purpose here has been to draw a more direct line from the geometry to Solèr's Theorem using wherever possible geometric methods intrinsically suited to the infinite-dimensional case. For a broader view of the subject as it existed prior to Solèr's Theorem, one should consult the unexcelled expert treatment given by S. Maeda.

**Definition 3.1** [V-Y, p. 16; Bir, Chapter IV, §7; M-M, §16]. A **projective geometry** $\mathbb{P}$ is a set, of cardinality at least 2, whose elements are called **points**, together with a distinguished family of subsets of $\mathbb{P}$, called **lines**, that together satisfy these three conditions:

**Axiom 1.** *For every two different points in $\mathbb{P}$ there is one, and only one, line that contains both these points.*

**Axiom 2.** *Every line contains at least three different points.*

The third axiom is best expressed using the following notation and terminology. If P and Q are different points in $\mathbb{P}$, we denote by $P+Q$ the unique line that contains them both. If $P = Q$, we set $P+Q = P = Q$. If the points $P$, $Q$, $R$ do not lie on the same line, then the three distinct lines $P+Q$, $Q+R$, $R+P$ constitute a **triangle**.

**Axiom 3.** *If a line intersects two sides of a triangle at different points, then it intersects the third side.* (*In more detail*: *If $P$, $Q$, $R$ are three noncollinear points and if the point $S$ lies on $P+Q$, the point $T$ on $Q+R$, and $S \neq T$, then $S+T$ and $R+P$ intersect.*)

Starting with these three axioms, one builds up the theory of **synthetic projective geometry**.

A subset $\ell \subseteq \mathbb{P}$ is a **subspace** if for every $P, Q \in \ell$ we have $P+Q \subseteq \ell$. By convention the empty set (denoted 0) and the points are also subspaces. Given two nonempty subspaces $\ell_1$, $\ell_2$ of $\mathbb{R}$, define

$$\ell_1 + \ell_2 = \{X \in \mathbb{P} : \ X \in P+Q \text{ for some } P \in \ell_1, Q \in \ell_2\}.$$

Thus $\ell_1 + \ell_2$ consists of all those points in $\mathbb{P}$ that lie on some line joining a point in $\ell_1$ with one in $\ell_2$. Also set $0 + \ell = \ell + 0 = \ell$ for every subspace $\ell$. It is a nontrivial fact that $\ell_1 + \ell_2$ is always a subspace, and that the associative law holds: $\ell_1 + (\ell_2 + \ell_3) = (\ell_1 + \ell_2) + \ell_3$ for all subspaces $\ell_1, \ell_2, \ell_3$ of $\mathbb{P}$ [M–M, Lemma 16.2]. This fact allows us to define the finite sum of subspaces, $\ell_1 + \ell_2 + \cdots + \ell_n$, and the infinite sum $\sum(\ell_\alpha : \ \alpha \in A)$ which is the set union of all finite subsums. In each case the result is a subspace.

A finite or infinite family $\{P_\alpha : \ \alpha \in A\}$ of two or more points is **independent** [Ba, VII.2; vN, Part I, Chapter II] if $\sum(P_\beta : \ \beta \in B) \cap \sum(P_\gamma : \ \gamma \in C) = 0$ for any decomposition $A = B \cup C$ of the indexing set $A$ into nonempty disjoint subsets $B$ and $C$. It is equivalent to require that no one of the points $P_\alpha$ lies in the subspace generated by the others. The family $\{P_\alpha : \ \alpha \in A\}$ is independent if, and only if, every finite subfamily is independent. Every nonzero subspace $\ell$ of the projective geometry $\mathbb{P}$, in particular $\ell = \mathbb{P}$ itself, has a **basis**, i.e. an independent family of points that generates $\ell$. All bases of a given subspace $\ell$ contain the same (cardinal) number of points [Ba, VII.2, Theorems 1 and 2]. This invariant is called the **rank** of $\ell$, and we add the convention that rank $(0) = 0$. Define the **dimension** of $\ell$, $\dim(\ell)$, as $\text{rank}(\ell) - 1$. Thus $\dim(0) = -1$, $\dim(point) = 0$, $\dim(line) = 1$. A subspace of dimension 2 we call a **plane**. (When it is necessary to distinguish this concept of rank or dimension from another, we shall specify this as **projective** **rank** or **projective dimension**.) For an infinite-dimensional subspace $\ell$ we shall write simply $\dim(\ell) = \infty$ when knowledge of the cardinal number is not needed. With this convention and the usual conventions about the symbol $\infty$, the formula $\dim(\ell + m) + \dim(l \cap m) = \dim(\ell) + \dim(m)$ holds generally for subspaces $\ell$, $m$.

Let $\{P_\alpha : \alpha \in A\}$ be a basis for the projective geometry $\mathbb{P}$. Given a point $X \in \mathbb{P}$, there are finitely many basis points, say, $P_1, P_2, \cdots, P_n$, such that $X \in P_1 + P_2 + \cdots + P_n$, and, if $n$ is minimal, that is, if $X$ belongs to no subspace generated by fewer than $n$ basis elements, then the $P_i$ occurring there are unique. Thus to every point $X$ in $\mathbb{P} = \sum(P_\alpha : \alpha \in A)$ we associate a unique finite subset of the basis $\{P_\alpha\}$ which generates a subspace containing $X$.

An **isomorphism** or **projectivity** $\varphi$ of the projective geometry $\mathbb{P}_1$ onto the projective geometry $\mathbb{P}_2$ is a one-to-one map of the point set $\mathbb{P}_1$ onto the point set $\mathbb{P}_2$ that satisfies $\varphi(P + Q) = \varphi(P) + \varphi(Q)$ for all points $P$, $Q$ in $\mathbb{P}_1$. One verifies easily that if $\varphi$ is an isomorphism of $\mathbb{P}_1$ onto $\mathbb{P}_2$, then $\varphi$ maps the family of all subspaces of $\mathbb{P}_1$ one-to-one onto the family of all subspaces of $\mathbb{P}_2$ and preserves inclusion: $\ell \subseteq m \iff \varphi(\ell) \subseteq \varphi(m)$. If $\varphi$ is an isomorphism, so is $\varphi^{-1}$.

If $m$ is a subspace of the projective geometry $\mathbb{P}$, $\dim(m) \geq 1$, then we denote by $[0, m]$ the set of points of $\mathbb{P}$ contained in $m$, together with those lines of $\mathbb{P}$ determined by pairs of different points in $m$. Each of these lines is a subset of $m$, and this collection of points and lines taken together satisfies our three axioms; hence $[0, m]$ is a projective geometry in its own right. The subspaces of $[0, m]$ are exactly those subspaces $\ell$ of $\mathbb{P}$ that satisfy $\ell \subseteq m$.

Given subspaces $\ell$, $m$ of $\mathbb{P}$, $\dim(\ell) \geq 1$, suppose that there exists a third subspace $x$ satisfying $\ell + x = m + x$, $\ell \cap x = m \cap x = 0$. Then the map $\varphi(P) = (P + x) \cap m$ is an isomorphism of the projective goemetry $[0, \ell]$ onto the projective geometry $[0, m]$. This map is called a **perspectivity with axis** $x$.

Given a subspace $\ell$ of the projective geometry $\mathbb{P}$, a subspace $x$ that satisfies $\ell + x = \mathbb{P}$, $\ell \cap x = 0$ is a **complement** of $\ell$. Every subspace has a complement and all complements of $\ell$ have the same dimension, because if $x_1$ and $x_2$ are complements of $\ell$, then $\ell$ serves as an axis of perspectivity between $[0, x_1]$ and $[0, x_2]$. The dimension of a complement of $\ell$ is called the **codimension** of $\ell$. A subspace of codimension $0$ is called a **hyperplane**. Thus $\ell$ is a hyperplane when there exists a point $P$ not on $\ell$ such that $\ell + P = \mathbb{P}$. The term "hyperplane" is thus a misnomer; "dual-point" or even "hyperpoint" would be more apt.

That is the basic material of synthetic projective geometry, whether finite-or infinite-dimensional. With this material in hand, one goes on to derive the various geometric facts that constitute the subject matter of synthetic projective geometry. For example, the theorem of Desargue can be put this way. Let $p$ and $p'$ be two different perspective intersecting subspaces of finite dimension $\geq 2$ in a projective geometry. Let $\ell$ be a subspace of $p$, $\ell' \subseteq p'$ its image under the given perspectivity. Then the assertion of Desargue's Theorem is embodied in the equation

$$\ell \cap \ell' = \ell \cap (p \cap p').$$

In the classical case, $p$ and $p'$ are different perspective planes that intersect in a line $p \cap p'$. Then the preceding displayed equation says that any line $\ell$ in $p$, $\ell \neq p \cap p'$ intersects its perspective image $\ell'$ (in $p'$) in the point $\ell \cap (p \cap p')$. Thus all such intersection points are collinear—they all lie on the line $p \cap p'$. Put this way, Desargue's Theorem has little to do with triangles.

The next definition, which generalizes the classical projective-geometric concept of polarity to the infinite-dimensional case, is in accord with a suggestion made by Professor Garrett Birkhoff thirty years ago (1964).

**Definition 3.2.** A **polarity** on a projective geometry $\mathbb{P}$ is a map $\ell \longmapsto \ell^\perp$ of the family of all subspaces of $\mathbb{P}$ into itself that satisfies the following three conditions:

1) $\mathbb{P}^\perp = 0$.
2) $\ell \subseteq m \Longrightarrow \ell^\perp \supseteq m^\perp$.
3) If $P$ is a point, then $P^\perp$ is a hyperplane, and $P^{\perp\perp} = P$.

We shall be dealing henceforth primarily with projective geometries that carry a polarity, often using the notation $\{\mathbb{P}, \perp\}$ to put explicit emphasis on this fact.

Given a projective geometry-with-polarity $\{\mathbb{P}, \perp\}$, it is natural to extend the classical terminology and call the hyperplane $P^\perp$ the **polar** of the point $P$; $P$ the **pole** of $P^\perp$. One establishes directly from (2) and (3) that $\ell \subseteq \ell^{\perp\perp}$ for every subspace $\ell$ and that $(\sum \ell_\alpha)^\perp = \cap \ell_\alpha^\perp$ for every family of subspaces $\{\ell_\alpha\}$. We call a subspace $\ell$ of $\{\mathbb{P}, \perp\}$ **closed** with respect to the given polarity when $\ell = \ell^{\perp\perp}$. Condition (3) states that every polar $P^\perp$ is closed; condition (3) also implies that every finite-dimensional subspace is closed. Hence, in the finite-dimensional case, our definition of polarity coincides with the classical definition: a bijective, order-inverting, period-2 map of the family of all subspaces of $\mathbb{P}$ onto itself [Ba, IV.3].

An infinite-dimensional projective geometry $\{\mathbb{P}, \perp\}$ always has nonclosed subspaces [Ba, IV.1, Existence Theorem]. The map $\ell \longmapsto \ell^{\perp\perp}$ is a kind of closure operation on the subspaces of $\mathbb{P}$. By way of motivation of the key definition to follow, let us consider how this closure operation can be relativized to a subgeometry $[0, m]$. Assuming that $m \cap m^\perp = 0$ (we call such a subspace **nonsingular**) and assuming $\dim(m) \geq 1$, one may check directly that the map $\ell \longmapsto \ell' = \ell^\perp \cap m$ is a polarity on the projective geometry $[0, m]$. Along with this polarity $\ell \longmapsto \ell'$ on $[0, m]$ comes its attendant closure operation $\ell \longmapsto \ell'' = (\ell')' = (\ell^\perp \cap m)' = (\ell^\perp \cap m)^\perp \cap m$. Will this be the same as $\ell^{\perp\perp} \cap m$? Only exceptionally. This requirement, suitably modified, forms the basis for our key definition. A desire for economy of notation suggests that we use the same name for the same thing.

**Lemma 3.3.** *If a projective geometry-with-polarity $\{\mathbb{P}, \perp\}$ satisfies any one of the following conditions, then it satisfies them all (the symbols $\ell$, $m$, $n$ denote subspace of $\mathbb{P}$):*

1) $\ell \subseteq m$, $m = m^{\perp\perp} \Longrightarrow (\ell^\perp \cap m)^\perp \cap m = \ell^{\perp\perp}$.
2) $\ell = \ell^{\perp\perp} \Longrightarrow \ell + \ell^\perp = \mathbb{P}$.
3) $m \subseteq n$, $m = m^{\perp\perp} \Longrightarrow m + (n \cap m^\perp) = n$.

**Definition 3.4.** We call $\{\mathbb{P}, \perp\}$ **orthomodular** when it satisfies the common condition of Lemma 3.3.

*Proof of Lemma* 3.3. (1) $\Longrightarrow$ (2) : If $\ell + \ell^\perp \neq \mathbb{P}$ for some closed subspace $\ell$, then there is a point $P \notin \ell + \ell^\perp$. Set $m = \ell + P \neq \ell$; $m$ is also closed. Claim $m \cap \ell^\perp = 0$, because otherwise there is a point $Q \in m \cap \ell^\perp$; so $Q \in m = \ell + P$, whence $Q \in L + P$ for some point $L \in \ell$. As $Q$ lies on the line $L + P$, then $P \in Q + L$, which puts $P$ in the subspace $\ell + \ell^\perp$, a contradiction. Hence, $m \cap \ell^\perp = 0$, so, by (1), $(\ell^\perp \cap m)^\perp \cap m = 0^\perp \cap m = \mathbb{P} \cap m = m = \ell^{\perp\perp} = \ell$, contradicting $m \neq \ell$. Hence $\ell + \ell^\perp = \mathbb{P}$ for every closed subspace $\ell$.

(2) $\Longrightarrow$ (1): Consider the left side of (1), remembering $m = m^{\perp\perp}$: $(\ell^\perp \cap m)^\perp \cap m = (\ell^\perp \cap m)^\perp \cap m^{\perp\perp} = ((\ell^\perp \cap m) + m^\perp)^\perp = ((\ell^\perp + m^\perp) \cap (m + m^\perp))^\perp = (\ell^\perp \cap \mathbb{P})^\perp = \ell^{\perp\perp}$.

To get (3) $\Longrightarrow$ (2), put $n = \mathbb{P}$ in (3). For (2) $\Longrightarrow$ (3), intersect both sides of $m + m^\perp = \mathbb{P}$ with $n$.

A closed subspace of an orthomodular projective geometry $\{\mathbb{P}, \perp\}$ remains closed with respect to any closed subspace containing it; i.e. $\ell = \ell^{\perp\perp} \subseteq m = m^{\perp\perp} \Longrightarrow \ell = \ell''$. Further, if $\{\mathbb{P}, \perp\}$ is orthomodular and $m$ is a closed subspace of $\mathbb{P}$ with $\dim(m) \geq 1$, then the geometry $[0, m]$ is also orthomodular with respect to the relativized polarity $\ell \longmapsto \ell' = \ell^\perp \cap m$ (for $\ell \subseteq m$). Finally, in an orthomodular projective geometry every subspace is nonsingular $(m \cap m^\perp = 0)$. In the finite-dimensional case, this condition is equivalent to orthomodularity. Thus in the finite-dimensional case we really have nothing new.

The universal example of a projective geometry is constructed as follows: Start with a field $\mathcal{K}$ (remember that "field" for us is a commutative or noncommutative division ring), and a left vector space $E$ over $\mathcal{K}$. We require only $\dim(E) \geq 2$, $\dim(E) = \infty$ is allowed. If we interpret the one-dimensional subspaces of $E$ as the points of a projective geometry $\mathbb{P}(E, \mathcal{K})$, interpret the 2-dimensional subspaces of $E$ as the lines of $\mathbb{P}(E, \mathcal{K})$, and use ordinary set inclusion, then our three axioms are satisfied. The subspaces of $E$ correspond exactly to the subspaces of $\mathbb{P}(E, \mathcal{K})$, vector space dimension equalling projective rank. Thus, while the vector space dimension of a plane in $E$ is 2, it corresponds to a line in $\mathbb{P}(E, \mathcal{K})$ of projective dimension 1.

The determinative "Fundamental Theorem of Projective Geometry" asserts the universality of this example.

**Theorem 3.5** (Existence). *Given a synthetic projective geometry $\mathbb{P}$ of projective dimension $\geq 3$ ($\infty$ allowed), there is a field $\mathcal{K}$ and a left vector space $E$ over $\mathcal{K}$ such that $\mathbb{P}$ is isomorphic to $\mathbb{P}(E, \mathcal{K})$.*

(Uniqueness). *Suppose that $E$ is a left vector space over a field $\mathcal{K}$, $3 \leq \dim(E) \leq \infty$, that $F$ is a left vector space over a field $\mathcal{L}$, and that $\varphi$ is a projective geometry isomorphism of $\mathbb{P}(E, \mathcal{K})$ onto $\mathbb{P}(F, \mathcal{L})$. Then there is a field isomorphism $g$ of $\mathcal{K}$ onto $\mathcal{L}$ and a bijective $g$-linear transformation $A$ of $E$ onto $F$ that implements $\varphi$ in the sense that $\varphi(\mathcal{K}x) = \mathcal{L}A(x)$ for all $x \in E$. Suppose further that $h$ is another field isomorphism of $\mathcal{K}$ onto $\mathcal{L}$ and that $B$ is a bijective $h$-linear transformation of $E$ onto $F$. Then $B$ also implements $\varphi \iff$ there exists $\lambda \in \mathcal{L}$, $\lambda \neq 0$ such that $B(x) = \lambda A(x)$ for all $x \in E$. In this case $h(\rho) = \lambda g(\rho)\lambda^{-1}$ for all $\rho \in \mathcal{K}$.*

(The additive transformation $A$ of $E$ onto $F$ is $g$-linear when $A(\rho x) = g(\rho)A(x)$ for all $x \in E$, $\rho \in \mathcal{K}$.)

The existence part of the fundamental theorem lies deep, and complete proofs are rare, especially in the infinite-dimensional case. Baer gives a complete and fully detailed proof in the last chapter of his book [Ba]. His nontraditional 47-page proof, couched in the language of lattice theory, has been elegantly summarized by S. Maeda [M–M, §33]. Von Neumann has generalized the fundamental theorem to complemented modular lattices [vN, Part II, Theorem 14.1]. Von Neumann's proof can be specialized to yield a more-or-less traditional coordinatization for an infinite-dimensional projective geometry; here is a very brief sketch of such a proof.

First, select a basis $\{P_\alpha \, ; \, \alpha \in A\}$ of $\mathbb{P}$ indexed by the well-ordered set $A = \{0, 1, \cdots\}$. Normalize this basis by selecting on each line $P_\alpha + P_\beta (\alpha \neq \beta)$ a third point $C_{\alpha\beta} = C_{\beta\alpha}$ different from $P_\alpha$ and $P_\beta$ such that

$$(3.1) \qquad (C_{\alpha\beta} + C_{\beta\gamma}) \cap (P_\alpha + P_\gamma) = C_{\alpha\gamma}$$

for every triple $\alpha, \beta, \gamma$ of distinct elements from the indexing set $A$. This can be done using von Neumann's result [vN, Part II, Lemma 5.3] together with transfinite induction.

Each $C_{\alpha\beta}$ is then an axis of perspectivity between any pair of lines $P_\alpha + P_\gamma$, $P_\beta + P_\gamma$ that pass through a common hinge point $P_\gamma$. Therefore the map $\theta_{\alpha\beta} : P_\alpha + P_\gamma \to P_\beta + P_\gamma$ given by $\theta_{\alpha\beta}(X) = (X + C_{\alpha\beta}) \cap (P_\beta + P_\gamma)$, $X \in P_\alpha + P_\gamma$, is an isomorphism.

Second, make the points of the line $P_0 + P_1 \setminus P_1$ into a field by defining addition and multiplication in this set by the classic projective-geometric constructions. This is done for example in the book of Veblen and Young [V-Y, Vol. 1, Chapter VI]. They state that their procedure is a clarification and simplification of the original algebra of throws of von Staudt. Of course the field $\mathcal{K}$ will not in general be commutative. Carry out this construction so that $P_0$ is the zero of $\mathcal{K}$, and $P_1$ is the point at infinity. At this stage, every point of $P_0 + P_1$ will have been assigned homogeneous coordinates $[\rho, \sigma, 0, \cdots]$ where $\rho, \sigma \in \mathcal{K}$ are not both zero. The place entry within the brackets correspond to the elements of the well-ordered indexing set $A$, and the square brackets $[\cdots]$ denote the equivalence class generated by left multiplication by nonzero elements of $\mathcal{K}$.

Third, using the isomorphisms $\theta_{\alpha\beta}$, transplant the homogeneous coordinates on $P_0 + P_1$ to every line $P_\alpha + P_\beta$ in a consistent way. The proof of consistency needs two identities satisfied by the maps $\theta_{\alpha\beta}$; these derive from von Neumann's normalization condition (3.1). At this point homogenous coordinates will have been introduced on every line $P_\alpha + P_\beta$. The general point on $P_\alpha + P_\beta$ has homogeneous coordinate $[\cdots 0, \rho, 0, \cdots, 0, \sigma, 0, \cdots]$ where the nonzero field elements $\rho, \sigma$ occur in the $\alpha$-place and $\beta$-place respectively and are not both zero. As every point $X \in \mathbb{P}$ lies in a unique subspace generated by finitely many basis points, the assignment of coordinates to $X$ becomes a standard construction of finite-dimensional projective geometry; see [V-Y, Vol.1, §§63, 70].

Finally, construct a left vector space $E$ over $\mathcal{K}$ with basis $\{e_\alpha \, ; \, \alpha \in A\}$ indexed by the same well-ordered set $A$. The one-dimensional subspaces $\mathcal{K}e_\alpha$ of $E$ are independent points in the projective geometry $\mathbb{P}(E, \mathcal{K})$. Every point in $\mathbb{P}(E, \mathcal{K})$ has the form $\mathcal{K}x$ for some nonzero vector $x \in E$. Express $x$ in terms of the basis: $x = \rho(1)e_{\alpha(1)} + \cdots + \rho(n)e_{\alpha(n)}$ with $0 \neq \rho(i) \in \mathcal{K}$. The projective point $\mathcal{K}x$ has homogenous coordinate $[\cdots, \rho(1), \cdots, \rho(2), \cdots, \text{etc.}]$ where $\rho(1)$ lies in the $\alpha(1)$-place, $\rho(2)$ in the $\alpha(2)$-place, etc.; zeroes elsewhere. Define the map $\varphi : \mathbb{P} \to \mathbb{P}(E, \mathcal{K})$ by assigning to $X \in \mathbb{P}$ the point $\mathcal{K}x$ in $\mathbb{P}(E, \mathcal{K})$ with the same homogenous coordinate. The proof that $\varphi(X + Y) = \varphi(X) + \varphi(Y)$ follows as in the finite-dimensional case.

The uniqueness part of the fundamental theorem is more straightforward, and many authors offer complete proofs; see, for example, [Ba, III.1; G-W, Chapter III, Theorem.6; or A, Theorem 2.26]. Note that uniqueness holds for projective dimension 2, while existence needs projective dimension $\geq 3$.

Thus every synthetic projective geometry $\mathbb{P}$ of projective dimension $\geq 3$, $\dim(\mathbb{P}) = \infty$ included, can be realized as a $\mathbb{P}(E, \mathcal{K})$, the family of all subspaces of a left vector space $E$ over an essentially unique field $\mathcal{K}$. In $\mathbb{P}(E, \mathcal{K})$ the projective-geometric operations $+, \cap$ become ordinary sum and intersection of subpsaces. If $\mathbb{P}$ carries a polarity $\perp$, then this polarity, transferred to $\mathbb{P}(E, \mathcal{K})$, becomes a map of the family of all subspaces of $E$ into itself. Such a map on the family of subspaces of the $\mathcal{K}$-vector space $E$ always arises as the orthogonality relation induced by a form on $E$—this is essentially the content of the famous theorem of Birkhoff and von Neumann.

**Theorem 3.6** [B-vN, §14 and Appendix; Ba, Chapter IV]. *Let $\mathcal{K}$ be a field of characteristic $\neq 2$, and let $E$ be a left vector space over $\mathcal{K}$ of dimension $\geq 3$, $\dim(E) = \infty$ included. Suppose given a polarity $\perp$ on $\mathbb{P}(E, \mathcal{K})$; that is, suppose given a map $L \longmapsto L^\perp$ of the family of all subspaces of $E$ into itself that satisfies the following three conditions*: (1) $E^\perp = 0$, (2) $L \subseteq M \Longrightarrow L^\perp \supseteq M^\perp$, *and* (3) $0 \neq x \in E \Rightarrow (\mathcal{K}x)^\perp$ *is a hyperplane in* $E$ *and* $(\mathcal{K}x)^{\perp\perp} = \mathcal{K}x$.

*Then either* (H): *There exists an* $e \in E$ *such that* $\mathcal{K}e \not\subseteq (\mathcal{K}e)^\perp$ *(the polarity is* **Hermitian***). In this case, $\mathcal{K}$ has an involution $\rho \longmapsto \rho^*$, and $E$ carries a $*$-Hermitian form $< \cdot, \cdot >$ whose orthogonality coincides with the given polarity. Moreover if $\rho \longmapsto \rho^\#$ is another involution on $\mathcal{K}$ and $[\cdot, \cdot]$ is a $\#$-Hermitian form that induces the same polarity, then there exists $0 \neq \lambda = \lambda^* \in \mathcal{K}$ such that $\rho^\# = \lambda^{-1}\rho^*\lambda$ for all $\rho \in \mathcal{K}$ and $[x, y] = < x, y > \lambda$ for all $x, y \in E$ and all $\rho \in \mathcal{K}$. $\mathbb{P}(E, \mathcal{K})$ is orthomodular (Definition 3.4) $\Longleftrightarrow \{E, \mathcal{K}, < \cdot, \cdot >\}$ is orthomodular (Definition 1.2).*

*or* (S): $\mathcal{K}x \subseteq (\mathcal{K}x)^\perp$ *for all* $x \in E$ *(the polarity is* **symplectic** *or* **alternate***). In this case $\mathcal{K}$ is commutative, and $E$ carries a symmetric form $< \cdot, \cdot >$ whose orthogonality coincides with the given polarity and which satisfies $< x, x > = 0$ for all $x \in E$. Moreover, if $[\cdot, \cdot]$ is another symmetric form on the commutative field $\mathcal{K}$ that induces the same polarity, then there exists $0 \neq \lambda \in \mathcal{K}$ such that $[x, y] = < x, y > \lambda$ for all $x, y$ in $E$.*

The kind of uniqueness described for the forms in cases (H) and (S) we can refer to as "unique up to scaling".

Birkhoff and von Neumann prove this theorem in the finite-dimensional case for a special kind of polarity called an orthocomplementation. Maeda gives an outline of their proof [M-M, Remark 34.3] and shows how to deduce the result in infinite dimensions from their finite-dimensional theorem [M–M, Theorem 34.5]. Baer studies general polarities on finite-dimensional spaces in full detail [Ba, Chapter IV]. As Theorem 3.6 is the really essential link in our application of Solér's Theorem to projective geometry, we sketch a proof that applies directly to the infinite-dimensional case.

The full dual $E'$ of $E$ consists of all $\mathcal{K}$-linear functionals $f : E \longrightarrow \mathcal{K}$    :

$f(\rho x) = \rho f(x)$ for all $x \in E$, $\rho \in \mathcal{K}$. The dual space $E'$ is a right vector space over $\mathcal{K}$, and a left vector space over $\mathcal{K}^{opp}$, the field $\mathcal{K}$ with multiplication reversed. Denoting multiplication in $\mathcal{K}^{opp}$ by a centered dot, $\rho \cdot \sigma = \sigma \rho$, we have $(\rho \cdot f)(x) = \rho \cdot f(x) = f(x)\rho$ for all $x \in E$, $f \in E'$, $\rho \in \mathcal{K}^{opp}$. Given a closed hyperplane $M = (\mathcal{K}x)^\perp$ in $E$, there is a functional $f_M \in E'$ that satisfies $\ker(f_M) = (\mathcal{K}x)^\perp$ and is uniquely determined up to a scalar multiple. The $f_M$ that arise this way generate a subspace $F$ of $E'$. Construct the projective geometry $\mathbb{P}(F, \mathcal{K}^{opp})$. Denote by $\varphi$ the map that takes the one-dimensional subspace $\mathcal{K}x$ of $E$ to the one-dimensional subspace $\mathcal{K}^{opp} \cdot f_M$ where $M = (\mathcal{K}x)^\perp$. One may check that $\varphi$ is a projective geometry isomorphism of $\mathbb{P}(E, \mathcal{K})$ onto $\mathbb{P}(F, \mathcal{K}^{opp})$.

As $\dim(E) \geq 3$, the uniqueness part of Theorem 3.5 applies, so that there is an isomorphism $g$ of $\mathcal{K}$ onto $\mathcal{K}^{opp}$ and a g-linear map $A$ of $E$ onto $F$ that implements $\varphi$ in the sense that $\varphi(\mathcal{K}x) = \mathcal{K}^{opp} \cdot A(x)$ for all $x \in E$. The map $g$, being an isomorphism of $\mathcal{K}$ onto $\mathcal{K}^{opp}$, is therefore an anti-automorphism of $\mathcal{K}$ onto itself. The map $A$ satisfies $A(\rho x) = g(\rho) \cdot A(x) = A(x)g(\rho)$ for all $x \in E$, $\rho \in \mathcal{K}$. Because the kernel of the functional $A(x)$ is $(\mathcal{K}x)^\perp$, we have $A(x)(y) = 0 \iff y \in (\mathcal{K}x)^\perp$.

Define $[x, y] = A(y)(x)$ for all $x, y$ in $E$. This $\mathcal{K}$-valued function is additive in each variable separately, linear in the first variable, and $g$-linear in the second ($[x, \rho y] = [x, y]g(\rho)$). Further, $[x, y] = 0 \iff [y, x] = 0$, because $\mathcal{K}x \subseteq (\mathcal{K}y)^\perp$ is equivalent to $\mathcal{K}y \subseteq (\mathcal{K}x)^\perp$. The orthogonality induced by this form coincides with the given polarity, because if $M$ is a subspace of $E$ with basis $\{e_\alpha : \alpha \in A\}$, then

$$M^\perp = \cap((\mathcal{K}e_\alpha)^\perp, \alpha \in A), \text{ so } x \in M^\perp \iff x \in (\mathcal{K}e_\alpha)^\perp \;\forall \alpha \iff$$
$$A(e_\alpha)(x) = 0 \;\forall \alpha \iff [x, e_\alpha] = 0 \;\forall \alpha \iff [x, m] = 0 \;\forall m \in M.$$

In particular, if $[a, x] = 0$ for all $x \in E$, then $a \in E^\perp = 0$ (Item 2 of Definition 1.1).

There are now two possibilities: either $[x, x] = 0$ for all $x \in E$, or not. In the first case it is a straighforward matter to show then that $g$=identity, $\mathcal{K}$ is commutative, and $[x, y] = -[y, x]$ for all $x, y$ in $E$. In this case set $< \cdot, \cdot >= [\cdot, \cdot]$ to establish case (S) of Theorem 3.6. In the second case, there is a vector $e \in E$ with $[e, e] \neq 0$. Set $\varepsilon = [e, e]$, define $< x, y >= [x, y]\varepsilon^{-1}$ for all $x, y$ in $E$, and set $\rho^* = \varepsilon g(\rho)\varepsilon^{-1}$ for all $\rho \in \mathcal{K}$. Then, arguing with some care, and using the fact that the characteristic of $\mathcal{K}$ is not 2, one concludes that the map $\rho \longmapsto \rho^*$ is an involution on $\mathcal{K}$, and that $< \cdot, \cdot >$ is Hermitian with respect to this involution [Definition 1.1(3)]. We have therefore proved that if case (S) of Theorem 3.6 does not obtain, then the Hermitian case (H) does. The fact that $\mathbb{P}(E, \mathcal{K})$ is orthomodular (Definition 3.4) exactly when $\{E, \mathcal{K}, < \cdot, \cdot >\}$ is orthomodular (Definition 1.2) follows from Lemma 3.3(2).

That completes the sketch of proof of the existence portion of Theorem 3.6; we shall omit the proof of uniqueness.

Theorems 3.5 and 3.6 taken together tell us that we may realize our infinite-dimensional synthetic orthomodular projective geometry-with-polarity $\{\mathbb{P}, \perp\}$ as $\mathbb{P}(E, \mathcal{K})$, the geometry of all subspace of a left vector space $E$ over a $*$-field $\mathcal{K}$, the polarity $\perp$ being that induced by an orthomodular form $< \cdot, \cdot >$ on $E$.

We wish now to apply Solèr's Theorem to characterize those geometries-with-polarity that correspond to real, complex, and quaternionic Hilbert space. Note that in the representation described in Theorem 3.6, the polarity is assumed given as part of the basic data; the representation is for the **combination**: projective geometry + polarity. To apply Solèr's Theorem, we need to characterize, in terms involving only the given polarity and projective-geometric constructions, both hypotheses of her theorem: (1) orthomodularity and (2) the existence of an infinite orthogonal sequence $\{e_i\}$ with $<e_i, e_i>=<e_j, e_j>$ for all $i, j$. As Definition 3.4 already presents orthomodularity in the required manner, it only remains to address hypothesis (2).

Given two nonzero vectors $e, f'$, if there exists $0 \neq \rho \in \mathcal{K}$ such that $<e, e>= \rho <f', f'> \rho^*$, then $e$ and $f = \rho f'$ satisfy $<e, e>=<f, f>$. The relation $\lambda \equiv \mu \iff \exists\, 0 \neq \rho \in \mathcal{K}$ such that $\lambda = \rho \mu \rho^*$ is an equivalence relation on the set of elements in $\mathcal{K}$ represented by the form. To secure our result, it is therefore enough to produce an orthogonal sequence of nonzero vectors $\{e_i : i = 1, 2, \cdots\}$ such that $<e_i, e_i> \equiv <e_j, e_j>$ for all $i, j$. This can be done using, along with the polarity, the projective-geometric notion of **harmonic conjugate**.

Let $P, Q$ be orthogonal points in our projective geometry $(P \in Q^\perp)$, and consider the line $P + Q$. Given a point $C$ on $P + Q$, $C \neq P, Q$, the harmonic conjugate of $C$ with respect to $P$ and $Q$ is a fourth point on $P + Q$ constructed as follows: Choose any point $X$ not on $P + Q$, choose $Y$ on the line $X + P$, $Y \neq X, P$, and construct the points $C' = (C + Y) \cap (Q + X)$ and $C'' = (P + C') \cap (Q + Y)$. Then $(X + C'') \cap (P + Q)$ is the harmonic conjugate of $C$ with respect to $P$ and $Q$ [V-Y, Vol.1, §31]. As char$(\mathcal{K}) \neq 2$ the harmonic conjugate is well defined and different from $C$. It does not depend on the choice of the construction points $X$ and $Y$.

Given that same point $C$ on $P + Q$, we may also construct the polar of $C$ relative to the geometry $[0, P + Q]$. This relative polar is $C^\perp \cap (P + Q)$, and it is also a point on $P + Q$ because the line $P + Q$ has projective dimension 1.

The projective points $P$ and $Q$ have the form $P = \mathcal{K}e$, $Q = \mathcal{K}f$ for nonzero orthogonal vectors $e, f$ in $E$. Every vector $\rho e$, $0 \neq \rho \in \mathcal{K}$, determines the same point $P$, and $<\rho e, \rho e> \equiv <e, e>$ for all nonzero $\rho$ in $\mathcal{K}$. Similarly for $Q$ and $f$. **In order that** $<e, e> \equiv <f, f>$ **it is necessary and sufficient that there exist a point** $C$ **on** $P + Q$ **whose harmonic conjugate coincides with its relative polar** $C^\perp \cap (P + Q)$, because every such $C$ is uniquely representable as $C = \mathcal{K}(e + \rho f)$ and its harmonic conjugate is the point $\mathcal{K}(e - \rho f)$, as a direct computation shows. If $\mathcal{K}(e - \rho f) = C^\perp \cap (P + Q)$, then $0 = <e - \rho f, e + \rho f> = <e, e> - \rho <f, f> \rho^*$, whence $<e, e> \equiv <f, f>$. We are now in a position to apply Solèr's Theorem.

**Theorem 3.7.** *If the synthetic projective geometry-with-polarity $\{\mathbb{P}, \perp\}$ is orthomodular and has a sequence of orthogonal points $\{P_i : i = 1, 2, \cdots\}$ such that every line $P_i + P_{i+1}$ contains a third point $C_i$ whose harmonic conjugate coincides with its relative polar, $i = 1, 2, \cdots$, then $\{\mathbb{P}, \perp\}$ is the family of all subspaces of a real, complex, or quaternionic Hilbert space $H$, the polarity $\perp$ being that given by the inner product in $H$.*

## 4. Orthomodular lattices

The orthomodular axiom first emerged in lattice theory. The family of projections in any von Neumann algebra constitute a complete orthocomplemented lattice, the orthocomplementation being given by $E^\perp = I - E$. Also, the simple identity $E \leq F \Longrightarrow F = E + (F - E)$ always holds in any such projection lattice. This identity, expressed in lattice-theoretic terms, becomes: $a \leq b \Longrightarrow b = a \vee (b \wedge a^\perp)$. In the late 1950's evidence began to mount that this lattice-theoretic axiom, by itself, could support a considerable body of theory. Thus emerged the class of abstract lattices called orthomodular, these being orthocomplemented lattices that satisfy what has come to be called the orthomodular axiom: $a \leq b \Longrightarrow b = a \vee (b \wedge a^\perp)$. And thus began the program of studying complete orthomodular lattices as the natural vehicle for generalizing projection lattices of von Neumann algebras. One can say that orthomodular lattices were designed to serve as algebraic versions of the projection lattices of von Neumann algebras much as Baer ∗-rings were intended to serve as algebraic models of the algebras themselves. For all these matters, see the 1966 survey article [H1].

Some properties valid in all projection lattices remain valid in the general orthomodular lattice, and many others do not. The question, What is "extra special" about the projection lattice of a von Neumann algebra as compared with the general orthomodular lattice? is the question form of what can be termed the **coordinatization problem:** Characterize, if possible, projection lattices among general orthomodular lattices. While this problem was implicit in [H1], the theory of orthomodular lattices was at that time (1966) really too primitive to justify proposing coordinatization as a feasible explicit project. But it arose as a central problem some fifteen years later.

Kalmbach, in her 1983 book *Orthomodular lattices* [Kal], which is a fully detailed comprehensive presentation of the subject up to that time, states as No. 29 of her "Problems on orthomodular lattices" the following

> *Characterize, lattice-theoretically, the orthomodular lattices of closed subspaces of infinite dimensional Hilbert spaces. In particular, consider these conditions on a complete orthomodular lattice: atomic, exchange axiom, infinite dimensional, irreducible. There are three examples, the lattice of closed subspaces of real, complex, quaternionic Hilbert space. Are there any other examples?* [Kal, p. 348]

As for that last question, Keller essentially provided the answer in 1980 [Ke]. The answer is yes, there are many other examples. Now, ten years later, the first problem is settled. We state the theorem first, then explain the terms involved.

**Theorem 4.1.** *Let $L$ be an irreducible, complete, orthomodular AC-lattice. If $L$ has an orthogonal sequence of atoms $\{p_i : i = 1, 2, \ldots\}$ together with another corresponding sequence of atoms $c_i \leq p_i \vee p_{i+1}$, $i = 1, 2, \ldots$ such that the harmonic conjugate of $c_i$ with respect to the pair of atoms $p_i$, $p_{i+1}$ equals $c_i^\perp \wedge (p_i \vee p_{i+1})$, $i = 1, 2, \ldots$, then $L$ is orthoisomorphic to the lattice of all closed subspaces of real, complex, or quaternionic Hilbert space.*

Now for the terms involved. An **orthomodular lattice** $L$ is, first of all, a lattice with 0 and 1, that is to say, a partially ordered set with smallest element 0 and largest element 1 in which every pair of elements $a$, $b$ has a least upper bound, or join, denoted $a \vee b$, and a greatest lower bound, or meet, denoted $a \wedge b$. Secondly $L$ carries an orthocomplementation which is a one-to-one map $a \mapsto a^\perp$ of $L$ onto itself that satisfies $a \wedge a^\perp = 0$, $a \vee a^\perp = 1$, $a \leq b \Longrightarrow a^\perp \geq b^\perp$, and $a^{\perp\perp} = a$ for all $a$, $b$ in $L$. And, finally, $L$ satisfies the orthomodular identity: $a \leq b \Longrightarrow b = a \vee (b \wedge a^\perp)$.

$L$ is **complete** if every nonempty set of elements of $L$ has a meet and join and is **irreducible** if its center consists of just 0 and 1. The center of $L$ consists of those $z \in L$ that commute with every $a \in L$, which means that $z = (z \wedge a) \vee (z \wedge a^\perp)$ for every $a \in L$. (Commutativity is symmetric; if $a$ commutes with $b$, then $b$ commutes with $a$.)

The term "AC" means **atomistic with the covering property** [M-M, Definition 8.7]. We say that $b$ covers $a$ when $b > a$ and $b > c > a$ is satisfied by no $c$. An element of $L$ that covers 0 is called an atom, and $L$ is atomistic when every $0 \neq a \in L$ is the join of all the atoms $p \leq a$. $L$ has the covering property when $a \wedge p = 0$ and $p$ an atom $\Longrightarrow a \vee p$ covers $a$.

Given orthogonal atoms $p$ and $q$ $(p \leq q^\perp)$, the **harmonic conjugate** of the atom $c \leq p \vee q$ is constructed by exactly the same procedure as in projective geometry: follow the procedure in the second paragraph preceding Theorem 3.7, replacing "point" by "atom", $+$ by $\vee$, and $\cap$ by $\wedge$.

That explains the meaning of all the terms in Theorem 4.1. Its proof rests on the characterization of real, complex, and quaternionic projective geometries (our Theorem 3.7). In the book [M-M], S. Maeda details, with clarity and precision, virtually all essential aspects of the proof, which in outline runs as follows:

1. Embed the lattice $L$ in a projective geometry $\mathbb{P}(L)$ whose points are exactly the atoms of $L$.
2. Use the given orthomodular orthocomplementation $a \longmapsto a^\perp$ on $L$ to induce an orthomodular polarity $\perp$ on $\mathbb{P}(L)$, observing that the elements of $L$ correspond precisely to the closed subspaces of $\{\mathbb{P}(L), \perp\}$.
3. Apply Theorem 3.7, which tells us that $\{\mathbb{P}(L), \perp\}$ is the family of all subspaces of a classical Hilbert space, whence $L$ appears as its lattice of closed subspaces.

Inasmuch as S. Maeda's proof is somewhat scattered throughout [M-M], which contains a wealth of other lattice-theoretic results, perhaps it will be a useful service to walk through the proof here, with specific references to the applicable results in [M-M] at each step.

1. The underlying set of points of our projective geometry $\mathbb{P}(L)$ consists of the atoms in $L$. The lines of $\mathbb{P}(L)$ are the subsets $\overline{ab} = \{p : p$ an atom in $L$ and $p \leq a \vee b\}$ determined by all pairs of different atoms $a$, $b$ in $L$ [M-M, Lemma 16.4]. To verify that $\mathbb{P}(L)$ is a projective geometry, we must check the three axioms of Definition 3.1: Axiom 1 follows from the covering property [M-M, combine Lemma 7.9 and 16.4 (I)]. Thus any pair of different points of $\mathbb{P}(L)$ lie on one, and only one, line. Axiom 2 is a consequence of irreducibility [M-M, combine 28.8.3 with 11.6; our lattice is DAC and finite-modular]. So every line in $\mathbb{P}(L)$ contains at least three different points. Axiom 3 is also a consequence of the covering property [M-M, Lemma 16.4 (II)]. Hence if a line in $\mathbb{P}(L)$ intersects two sides of a triangle at different points, then it intersects the third side. Thus $\mathbb{P}(L)$ is a projective geometry, and the entire machinery of §3 now comes into play: subspaces, sums and intersections of subspaces (both finite and infinite), independence of points, bases, and projective dimension. We use the notation of §3— $\in$, $\subseteq$, $+$, and $\bigcap$ in $\mathbb{P}(L)$—and continue to use $\leq$, $\vee$, and $\wedge$ for the operations in our lattice $L$. Given a nonzero element $a$ in $L$, we associate to it a subspace $\omega(a)$ of $\mathbb{P}(L)$ given by $\omega(a) = \{p \,:\, p \text{ an atom in } L \text{ and } p \leq a\}$. Also set $\omega(0) = 0$. The map $a \longmapsto \omega(a)$ is a one-to-one order-preserving mapping of our lattice $L$ into the projective geometry $\mathbb{P}(L)$ that satisfies $\omega(\bigwedge a_\alpha) = \bigcap \omega(a_\alpha)$ for any family $\{a_\alpha : \alpha \in A\}$ of elements from $L$ [M-M, Theorem 15.5]. The set of those elements of $L$ which are the join of finitely many atoms is closed under the operations $\vee$ and $\wedge$, and, for such elements, $\omega(a \vee b) = \omega(a) + \omega(b)$ [M-M, Definition 8.8 and Theorem 15.5(3)]. If $a^\perp$ is an atom in $L$, then $\omega(a)$ is a hyperplane in $\mathbb{P}(L)$ [M-M, Definition 7.1 and Lemma 17.17]. In this way we embed our lattice $L$ in the projective geometry $\mathbb{P}(L)$.

2. Given a point $p$ in $\mathbb{P}(L)$ ($p$ is an atom in $L$), define its polar to be the hyperplane $\omega(p^\perp) = \{r \,:\, r \text{ an atom and } r \leq p^\perp\}$ in $\mathbb{P}(L)$. Extend the polarity to all subspaces $\ell$ by $\ell^\perp = \bigcap(\omega(p^\perp) \,:\, p \in \ell)$. Then the double polar of $p$ is

$$\begin{aligned}(\omega(p^\perp))^\perp &= \bigcap (\omega(r^\perp) \,:\, r \text{ an atom and } r \in \omega(p^\perp)) \\ &= \bigcap (\omega(r^\perp) \,:\, r \text{ an atom and } r \leq p^\perp) \\ &= \omega (\bigwedge(r^\perp \,:\, r \text{ an atom and } r \leq p^\perp)) \\ &= \omega (\bigvee(r \,:\, r \text{ an atom and } r \leq p^\perp))^\perp \\ &= \omega (p^{\perp\perp}) = \omega(p) = p.\end{aligned}$$

Hence condition (3) of Definition 3.2 is satisfied. Conditions (1) and (2), namely, $\mathbb{P}(L)^\perp = 0$ and $\ell \subseteq m \Longrightarrow \ell^\perp \supseteq m^\perp$, are easily checked, so we have verified that the orthocomplementation on $L$ induces in a natural way a polarity $\perp$ on the associated projective geometry $\mathbb{P}(L)$. Next, observe that the closed subspaces of $\{\mathbb{P}(L), \perp\}$, namely, those that satisfy $\ell = \ell^{\perp\perp}$, are exactly those of the form $\ell = \omega(a)$ for some $a$ in $L$, because if $\ell$ is any subspace of $\mathbb{P}(L)$, then

$$\begin{aligned}\ell^\perp &= \bigcap (\omega(p^\perp) \,:\, p \text{ an atom and } p \in \ell) \\ &= \omega (\bigwedge(p^\perp \,:\, p \text{ an atom and } p \in \ell)) \\ &= \omega ((\bigvee(p \,:\, p \text{ an atom and } p \in \ell))^\perp) \\ &= \omega (a^\perp) \text{ where } a = \bigvee(p \,:\, p \in \ell).\end{aligned}$$

As the closed subspaces of $\{\mathbb{P}(L), \perp\}$ are exactly those of the form $\ell^\perp$, every closed subspace equals $\omega(a)$ for some $a$ in $L$. Conversely, if $\ell = \omega(a)$, then as a direct calculation shows, $\ell^\perp = \omega(a^\perp)$ and $\ell^{\perp\perp} = \omega(a^{\perp\perp}) = \omega(a) = \ell$. Hence under the embedding $a \longmapsto \omega(a)$ of our lattice $L$ into the projective geometry-with-polarity $\{\mathbb{P}(L), \perp\}$, the elements of $L$ correspond exactly to the closed subspace of $\{\mathbb{P}(L), \perp\}$.

Finally, we need to prove that the orthomodularity of the lattice $L$, that is, the condition $a \leq b \Rightarrow b = a \vee (b \wedge a^\perp)$, implies the orthomodularity of the projective geometry $\{\mathbb{P}(L), \perp\}$ (Definition 3.4). Setting $a = d^\perp$, $b = c^\perp$ in the orthomodular axiom in $L$, then taking orthocomplements, we convert it to the equivalent form $c \leq d \implies (c^\perp \wedge d)^\perp \wedge d = c$. Consider the condition (1) of Lemma 3.3: $\ell \subseteq m$ and $m = m^{\perp\perp} \implies (\ell^\perp \cap m)^\perp \cap m = \ell^{\perp\perp}$. As established in the preceding paragraph, we have $\ell^\perp = \omega(c^\perp)$, $\ell^{\perp\perp} = \omega(c)$, and $m = \omega(d)$. As $\ell^{\perp\perp} \subseteq m^{\perp\perp} = m$, we have $c \leq d$ so $(c^\perp \wedge d)^\perp \wedge d = c$. Apply $\omega$: $\omega[(c^\perp \wedge d)^\perp \wedge d] = (\omega(c^\perp \wedge d))^\perp \cap \omega(d) = (\omega(c^\perp) \cap \omega(d))^\perp \cap \omega(d) = (\ell^\perp \cap m)^\perp \cap m = \omega(c) = \ell^{\perp\perp}$. Thus $\{\mathbb{P}(L), \perp\}$ is orthomodular.

3. The condition on harmonic conjugates stated in Theorem 4.1 transfers to the same condition in $\{\mathbb{P}(L), \perp\}$. Then, by Theorem 3.7, the projective geometry $\{\mathbb{P}(L), \perp\}$ is the lattice of all subspaces of a real, complex, or quatermonic Hilbert space, and our lattice $L$ appears as its closed subspaces.

This completes our summary of proof for Theorem 4.1.

A noteworthy perspective gained from this discussion, which is exposed in full detail by S. Maeda in [M-M], is this: irreducible, complete, orthomodular AC-lattices and orthomodular projective geometries are really one and the same thing. Given an orthomodular projective geometry, its lattice of closed subspaces has all the properties listed. And, given an irreducible, complete, orthomodular AC-lattice, we have just shown that it is the lattice of closed subspaces of an orthomodular projective geometry.

Finally, we note that in place of the condition in Theorem 4.1 involving harmonic conjugates we could also use the "angle bisecting" axiom of Morash [Mo].

## 5. Mackey's quantum logic

In his book *Mathematical foundations of quantum mechanics* [M], Professor George W. Mackey proposes an axiom system for nonrelativistic quantum mechanics. His basic elements are two abstract sets $\mathcal{O}$ and $\mathcal{S}$ and the family of Borel subsets of the real numbers. The axioms are expressed in terms of a postulated function $p(A, \alpha, X)$ that assigns a real number $s$, $0 \leq s \leq 1$, to each triple $A \in \mathcal{O}$, $\alpha \in \mathcal{S}$ and Borel set $X$. The elements $A \in \mathcal{O}$ are to be thought of as observables, the $\alpha \in \mathcal{S}$ as states, and the real number $p(A, \alpha, X)$, $0 \leq p(A, \alpha, X) \leq 1$, as the probability that a measurement of the observable $A$ on a system in state $\alpha$ will yield a value in the Borel set $X$. The first six axioms run as follows:

**Axiom I** asserts that for any observable $A$ and any state $\alpha$, $p(A, \alpha, \emptyset) = 0$, $p(A, \alpha, \mathbb{R}) = 1$, and $p(A, \alpha, \cdot)$ is countably additive on the Borel sets. Thus, with $A$ and $\alpha$ fixed, $p(A, \alpha, \cdot)$ is a probablity measure on the Borel sets.

**Axiom II** asserts that if $p(A, \alpha, X) = p(B, \alpha, X)$ for all states $\alpha$ and all Borel sets $X$, then $A = B$. Likewise, if $p(A, \alpha, X) = p(A, \beta, X)$ for all observables $A$ and all $X$, then $\alpha = \beta$.

**Axiom III** says this: Given any observable $A$ and any real-valued Borel measurable function $f$, then there is another observable $B$ such that $p(B, \alpha, X) = p(A, \alpha, f^{-1}(X))$ for every state $\alpha$ and every Borel set $X$. Axiom II shows that $B$ is uniquely determined by $A$; set $B = f(A)$. Thus, given an observable $A$, this axiom provides for the construction of other observables such as $A^2$, $\sin(A)$, etc. Axiom III thus postulates that the set $\mathcal{O}$ of observables is closed under a functional calculus based on Borel functions. It has a reasonable physical interpretation: the observable $A^2$ (for example) has a measured value $\lambda^2$ whenever $A$ has the measured value $\lambda$. Consider the particular function $f(x) = 1$ for all $x \in \mathbb{R}$. The observable $I = f(A)$ is uniquely determined by the property that $p(I, \alpha, X) = p(A, \alpha, f^{-1}(X))$ for every state $\alpha$ and every Borel set $X$. But $f^{-1}(X) = \mathbb{R}$ if $1 \in X$, $f^{-1}(X) = \emptyset$ if $1 \notin X$, so $p(I, \alpha, X) = 1$ when $1 \in X$ and is $0$ otherwise. Thus the observable $I$ has this property: whenever a measurement of $I$ is made on the system in any state $\alpha$, the probability that the measured value of $I$ will lie in the Borel set $X$ is $1$ when $1 \in X$, and $0$ otherwise. Hence the probability $p(I, \alpha, X)$ equals either $0$ or $1$ and depends only on $X$, not on $\alpha$ or $A$. Whenever the observable $I$ is measured, the measured value will be $1$. In the same way we may construct the observable $O$ using the function $f(x) \equiv 0$. Measurement of $O$ always yields the value $0$. We have $I^2 = I$ and $O^2 = O$.

**Axiom IV** provides for the formation of finite or countably infinite convex linear combinations of states. Specifically, if $\alpha_1, \alpha_2, \ldots$ is a finite or countably infinite set of members of $\mathcal{S}$ and $t_i$ a corresponding set of positive real numbers of sum $1$, then Axiom IV postulates the existence of a state $\alpha \in \mathcal{S}$ such that $p(A, \alpha, X) = \sum t_i p(A, \alpha_i, X)$ for all observables $A \in \mathcal{O}$ and all Borel sets $X$. This state $\alpha$ is uniquely determined (Axiom II), and we write $\alpha = \sum t_i \alpha_i$. A state is *pure* if it is not a convex linear combination of two states different from itself. Axiom IV does not postulate the existence of pure states; thus, in particular, it does not require that the pure states generate $\mathcal{S}$ under convex linear combinations, which is the case in basic quantum mechanics.

To make this rather abstract process more tangible, let us examine how all this plays out in a familiar simple example—the linear oscillator [H5, §5.5].

The linear oscillator consists of a particle of mass $m$ which, moving in a straight line, is attracted toward a fixed point on the line by a force proportional to its distance $x$ from that point. Taking this line as the real axis, the configuration space for this system is $\mathbb{R}$, and its "state space" is $L^2(\mathbb{R})$. The observable quantities of the oscillator, the set $\mathcal{O}$, correspond to selfadjoint operators on $L^2(\mathbb{R})$. The total energy $E$ of the oscillator is an observable and corresponds to the selfadjoint operator $H$ on $L^2(\mathbb{R})$ determined by the differential expression $-(\hbar^2/2m)d^2/dx^2 + 2\pi^2 \nu_0^2 m x^2$ where $\nu_0$ is the natural frequency of the system (the differential expression of Schrödinger's equation). Functions $f(x) \in L^2(\mathbb{R})$ that satisfy $<f, f> = \int_{-\infty}^{\infty} |f|^2 dx = 1$ correspond to pure states of our oscillator. The set $\mathcal{S}$ includes the pure states and other states as well. But we shall concentrate here on the pure states. The selfadjoint operator $H$ has eigenvalues $E_n = (2n+1)E_0$, $n = 0, 1, 2, \ldots$, where $E_0 = (\frac{1}{2})h\nu_0$ is the zero-point energy. Measurement of the total energy $E$ of the oscillator will yield only one of the values $E_n$; no other measured values of the total energy are possible. Each of these eigenvalues $E_n$ is simple, and the corresponding normalized eigenfunctions $\Psi_n(x)$, $n = 0, 1, 2, \ldots$ constitute an orthonormal basis of $L^2(\mathbb{R})$. These eigenfunctions correspond to particular pure states of the oscillator, its eigenstates. In this example, the function $p(E, \Psi, X)$ represents the probability that a measurement of the total energy $E$ of the oscillator in the pure state $\Psi$ will yield a value in the Borel set $X$. If $X$ does not contain one of the eigenval-ues $E_n$, that probability will be zero. In general, $p(E, \Psi, X)$ is obtained this way: The spectral family of $H$ consists of the projections $P(\lambda)$, $-\infty < \lambda < \infty$, where $P(\lambda)$ is the projection on the closed subspace of $L^2(\mathbb{R})$ generated by $\Psi_0, \Psi_1, ..., \Psi_n$, $n \leq \lambda < n+1$. If the oscillator is in the pure state $\Psi$, the probability that a measured value of its total energy lies in the interval $X = [\lambda, \mu]$ is $p(E, \Psi, X) = <(P(\mu) - P(\lambda))(\Psi), \Psi>$. If we use the more sophisticated "projection-valued measure" form of the spectral theorem for our energy operator $H$, which assigns to each Borel set $X$ a projection $P(X)$ on $L^2(\mathbb{R})$ [B-C, Appendix B,2], then $p(E, \Psi, X) = <P(X)(\Psi), \Psi>$ for any Borel set $X$.

Return now to the general formulation. Mackey now singles out special types of observables he calls "questions". An observable $A$ is a **question** when, for every state $\alpha \in \mathcal{S}$, its associated probability measure $p(A, \alpha, \cdot)$ equals 1 on the two-element Borel set $\{0, 1\}$. Use the letter $Q$ to denote this special type of observable. Thus, if $Q$ is a question, then whatever the state $\alpha$, $p(Q, \alpha, X) = 1$ for every Borel set $X$ that contains both 0 and 1. Which is to say a measurement of a question $Q$ on a system in any state $\alpha$ will yield a value in $X$ with certainty (probability 1) whenever $X$ contains both 0 and 1. Likewise $p(Q, \alpha, X) = 0$ if $X$ contains neither 0 nor 1. If $X$ contains 1 but not 0, then $p(Q, \alpha, X) = s$, $0 \leq s \leq 1$ ($s$ depending on $\alpha$), and if $X$ contains 0 but not 1, $1-s$. The actual measured value of the question $Q$ will be either 1 or 0. For a system in state $\alpha$, the probability that the measured value will be 1 is $p(Q, \alpha, \{1\}) = s$; the result of measurement will be 0 with probability $1-s$. We shall denote the set of questions by $\mathcal{L}$; it is a subset of the set $\mathcal{O}$ of observables: $\mathcal{L} \subseteq \mathcal{O}$.

Given an observable $A \in \mathcal{O}$, the functional calculus of Axiom III permits the construction of a family of questions associated to $A$. If $F \subseteq \mathbb{R}$ is any Borel set and $f$ its characteristic function ($f(x) = 1$ when $x \in F$, $f(x) = 0$ when $x \notin F$), then $f(A)$ is a question. (Because the observable $Q = f(A)$ is determined by the condition $p(Q, \alpha, X) = p(A, \alpha, f^{-1}(X))$ for all Borel sets $X$, and if $X$ is the two-element set $\{0, 1\}$, then $f^{-1}(X) = \mathbb{R}$, so $p = 1$ on $\{0, 1\}$ as required.) This question $Q = f(A)$ depends on the observable $A$, and on the Borel set $F$, so might better be written $Q = Q(A, F)$. Earlier we constructed two special questions, $I = Q(A, \mathbb{R})$ and $O = Q(A, \emptyset)$, both independent of the observable $A$. In general, the measured value of the question $Q(A, F)$, as of any question, is either 1 or 0. For a system in a particular state $\alpha$, the probability that a measurement of $Q(A, F)$ will result in the value 1 equals the probability that a measurement of the observable $A$ will be in $F$. If $A$ and $B$ are two observables, and $Q(A, F) = Q(B, F)$ for all Borel sets $F$, then $A = B$ (by Axiom II). So the family of all questions $Q(A, F)$, $F \in$ Borel sets, associated to a given observable $A$ uniquely determines $A$. Moreover, every question $Q_1$ arises as $Q(A, F)$ for some observable $A$ and some Borel set $F$ in a trivial way. Given the question $Q_1$, we note that it is, by definition, an observable, so we may form the question $Q = Q(Q_1, \{1\})$. Then $Q_1 = Q$.

The focus now shifts from the set $\mathcal{O}$ of all observables to its subset $\mathcal{L}$, the questions. And the states are given a new interpretation.

For each state $\alpha$ define a function $m_\alpha : \mathcal{L} \longrightarrow \{s : 0 \leq s \leq 1\}$ by $m_\alpha(Q) = p(Q, \alpha, \{1\})$. Thus $m_\alpha(Q)$ is the probability that a measurement of the question $Q$ on the system in state $\alpha$ yields the value 1. From Axiom II follows these two facts: if $m_\alpha(Q_1) = m_\alpha(Q_2)$ for all states $\alpha$, then $Q_1 = Q_2$; and if $m_\alpha(Q) =$ $m_\beta(Q)$ for all questions $Q$, then $\alpha = \beta$. Hence the map $\alpha \longmapsto m_\alpha$ is one-to-one, and we may identify the set $\mathcal{S}$ of states with the family of functions $m_\alpha$. We shall denote this set of functions also by $\mathcal{S}$. So, at this juncture, we have singled out a special family $\mathcal{L}$ of observables, the questions, and have identified the states $\mathcal{S}$ with certain [0,1]-valued functions on the observables. Attention now focuses on the system $\{\mathcal{L}, \mathcal{S}\}$.

The functions $m_\alpha : \mathcal{L} \longrightarrow \{s : 0 \leq s \leq 1\}$ are used to give $\mathcal{L}$ the structure of an orthomodular partially ordered set [M, p. 64]. Axioms V and VI are then added [M, pp. 65–66]. The result is a system $\{\mathcal{L}, \mathcal{S}\}$ of questions and states, where $\mathcal{L}$ is a countably orthocomplete orthomodular partially ordered set and $\mathcal{S}$ is a full, strongly convex, family of probability measures on $\mathcal{L}$. In detail:

**Definition 5.1.** $\mathcal{L}$ **is a countably orthocomplete orthomodular partially ordered set** if (1) $\mathcal{L}$ is a partially ordered set with smallest element O and largest element I; (2) $\mathcal{L}$ carries a bijective map $a \longmapsto a^\perp$ that satisfies $a^{\perp\perp} = a$, and $a \leq b \Longrightarrow a^\perp \geq b^\perp$ for all $a, b \in \mathcal{L}$; (3) for every $a \in \mathcal{L}$ the join $a \vee a^\perp = I$ and the meet $a \wedge a^\perp = O$ both exist and have the value indicated; (4) given any sequence $a_i$, $i = 1, 2, \ldots$ of elements from $\mathcal{L}$ such that $a_i \leq a_j^\perp$ when $i \neq j$, the join $\bigvee a_i$ exists in $\mathcal{L}$; and (5) $\mathcal{L}$ is orthomodular: $a \leq b \Longrightarrow b = a \vee (b \wedge a^\perp)$.

We say that $a, b \in \mathcal{L}$ are **orthogonal** written $a \perp b$, when $a \leq b^\perp$. The family $\{a_i\}$ is orthogonal when $i \neq j \Longrightarrow a_i \perp b_j$.

**Definition 5.2.** $\mathcal{S}$ **is a full, strongly convex family of probability measures on** $\mathcal{L}$ if (1) each $m \in \mathcal{S}$ is a probability measure on $\mathcal{L}$, that is, $m : \mathcal{L} \longrightarrow \{s : 0 \leq s \leq 1\}$, $m(O) = 0$, $m(I) = 1$, and $m(\bigvee a_i) = \sum m(a_i)$ for any orthogonal family $\{a_i : i = 1, 2, \ldots\}$ of elements of $\mathcal{L}$; (2) $m(a) \leq m(b)$ for all $m \in \mathcal{S} \Longrightarrow a \leq b$ ("full"); and (3) $m_i \in \mathcal{S}$, $0 < t_i \in \mathbb{R}$, $i = 1, 2, \ldots$, and $\sum t_i = 1$ together imply $\sum t_i m_i \in \mathcal{S}$ ("strongly convex").

Mackey states that the notion of a system $\{\mathcal{L}, \mathcal{S}\}$ satisfying Definitions 5.1 and 5.2 is equivalent to the notion of a system $\{\mathcal{O}, \mathcal{S}, p\}$ satisfying Axioms I through VI [M, Theorem on p. 68]; Beltrametti and Cassinelli have provided a complete proof of that equivalence [B-C, §13.7 and references there]. Thus the first six axioms of this approach to the mathematical foundations of quantum mechanics can be summed up in a pair $\{\mathcal{L}, \mathcal{S}\}$ having the properties listed in Definitions 5.1 and 5.2. $\mathcal{L}$ is named the **logic** of the physical system [M, p. 68], with reference in this connection to the famous paper "The Logic of Quantum Mechanics" by Garrett Birkhoff and John von Neumann [B-vN].

**Axiom VII** appears at this juncture. This axiom, together with some comments that precede and follow it, read thus [M, pp. 71–72]:

> …Almost all modern quantum mechanics is based implicitly or explicitly on the following assumption which we shall state as an axiom:
> **Axiom** VII: *The partially ordered set of all questions in quantum mechanics is isomorphic to the partially ordered set of all closed subspaces of a separable, infinite dimensional Hilbert space.*

> *This axiom has a rather different character from Axioms* I *through* VI. *These all had some degree of physical naturalness and plausibility. Axiom* VII *seems entirely ad hoc.... Can we justify making it?... Ideally one would like to have a list of physically plausible assumptions from which one could deduce Axiom* VII.

This is the project we undertake here: to devise "...*a list of physically plausible assumptions from which one could deduce Axiom* VII."

In 1961, Zierler proposed such a list [Z]. His axioms run as follows: (numbering as in the original paper).

> (E4),(E5),(A),(ND): $\mathcal{L}$ is a separable atomic lattice, center($\mathcal{L}$) $\neq$ $\mathcal{L}$, and $I \in \mathcal{L}$ is not finite.
> (M),(H): If $a \in \mathcal{L}$ is finite, then $\mathcal{L}(0, a)$ is modular; if $a$, $b$ are finite elements of the same dimension, then $\mathcal{L}(0, a)$ and $\mathcal{L}(0, b)$ are isomorphic.
> (S2): If $0 \neq a \in \mathcal{L}$, then there exists $m \in \mathcal{S}$ with $m(a) = 1$.
> (S3): $m(a) = 0$ and $m(b) = 0$ together imply $m(a \vee b) = 0$.
> (C'), (C): For every finite $a \in \mathcal{L}$ and for each $i$, $0 \leq i \leq \dim(a)$, the set of elements $\{x \in \mathcal{L} : x \leq a \text{ and } \dim x = i\}$ is compact in the topology provided by the metric
> $$f(x, y) = \sup\{|m(x) - m(y)| \, : \, m \in \mathcal{S}\}.$$
> For each $i = 0, 1, 2, ...$ the set of finite elements in $\mathcal{L}$ of dimension $i$ is complete with respect to the same metric.
> (Co): For some finite $b$ and real interval I there exists a nonconstant function from I to $\mathcal{L}(0, b)$.

We shall use a few of Zieler's axioms. But, having Solèr's Theorem, we shall not need (C'), (C), and (Co).

Our first axiom combines Zierler's (E4) with a poset-substitute for (S3) [Z, pp. 1155–1556].

**Axiom A.** (1) $\mathcal{L}$ is separable: i.e. any orthogonal family of nonzero elements in $\mathcal{L}$ is at most countable. (2) If $m(a) = m(b) = 0$ for some $a, b \in \mathcal{L}$ and an $m \in \mathcal{S}$, then there exists $c \in \mathcal{L}$, $c \geq a$ and $c \geq b$ with $m(c) = 0$.

Axiom A2, in the form given here, is taken from [B–C, §11.4]; we must use this form rather than Zierler's (S3), as we are not assuming that $\mathcal{L}$ is a lattice. Zierler provides no physical basis for his (S3); neither can I provide any for A2. It would certainly be desirable to do so; either that, or prove that A2 is nonessential. The separability, A1, is noncontroversial.

From Axiom A we deduce

**Lemma 5.3.** *For each state $m \in \mathcal{S}$ there is a unique element $a \in \mathcal{L}$ satisfying $m(x) = 1 \iff x \geq a$.*

The proof of this well-known result is easily given [B-C, p. 298]. The element $a \in L$ specified in Lemma 5.3 is called the **support** of the state $m$. Thus the support of $m \in \mathcal{S}$ can be characterized as the smallest element $a \in \mathcal{L}$ such that $m(a) = 1$. The element $b = a^\perp$ is called the maximal null element for $m$. It is characterized by $m(x) = 0 \iff x \leq b$. One may check that if $m_1$, $m_2$, ... is a finite or countably infinite sequence of states with respective supports $b_1$, $b_2$, ... in $\mathcal{L}$, and if $b_i \perp b_j$ when $i \neq j$, then $b = \bigvee b_i$ is the support of $m = \sum t_i m_i$, where the $t_i$ are any positive real numbers of sum 1. Also, if $b_1 = $ support $(m_1)$ and $b_2 = $ support $(m_2)$, then $b_1 \vee b_2$ exists in $\mathcal{L}$ and equals the support of $m = t_1 m_1 + t_2 m_2$ where $t_1$ and $t_2$ are any two positive numbers of sum 1.

Our next two axioms deal with pure states and postulate the validity of the principle of superpositions of pure states. Recall that a state $m \in \mathcal{S}$ is **pure** if it cannot be written as a convex linear combination of two states different from itself. That is, $m$ is pure if there do not exist $m_1$, $m_2 \in \mathcal{S}$ different from $m$ such that $m = t_1 m_1 + t_2 m_2$ for positive real numbers $t_1$, $t_2$ of sum 1. Why do we bring in pure states and their principle of superposition? The reason is that a system $\{\mathcal{L}, \mathcal{S}\}$ as described in Definitions 5.1 and 5.2 applies equally to classical mechanics and quantum mechanics. So our list of axioms must, among other things, distinguish quantum mechanics from classical mechanics. What differentiates quantum mechanics from classical mechanics? I contend that the concept of pure state, together with the principle of superposition of pure states, captures the essence of quantum mechanics. The following paragraphs contain my arguments supporting that contention.

In the early years of quantum mechanics the only states considered were pure, and they were called simply "states". For the working physicist of today, pure states are the basic building blocks. Dirac, writing in the fourth (1957) edition of his classic book *The principles of quantum mechanics* [D] uses the term "state" as synonymous with "pure state". Von Neumann, in his *Mathematical foundations of quantum mechanics* [vN2], also uses the term "state" for what are now called pure states. In Chapter IV of his book he says this (pp. 295–296):

> *... But the statistical character may become even more prominent, if we do not even know what state is actually present—for example, when several states $\phi_1$, $\phi_2$, ... with the respective probabilities $w_1$, $w_2$, ... ($w_1 \geq 0$, $w_2 \geq 0$, ..., $w_1 + w_2 + ... = 1$) constitute the description...*

And goes on to say (pp. 296–297):

> *We note that we shall have to pay attention to these mixtures of states also, in addition to the individual states themselves...*

When von Neumann says "state" here, he means what is today called a pure state. What he refers to as "mixtures of states" have become today's states. Thus both Dirac and von Neumann treat pure states as fundamental. They are Dirac's exclusive concern. And von Neumann clearly sees them as fundamental, with "mixtures of states" as a derived concept.

The distinction between pure states and states is nicely seen in our example of the linear oscillator. We have said that a given pure state of the oscillator corresponds to a function $\Psi \in L^2(\mathbb{R})$ of length $1 : <\Psi, \Psi> = \int_{-\infty}^{\infty} |\Psi|^2 dx = 1$. Actually, all functions obtained from $\Psi$ by multiplication by a complex number of absolute value 1 also correspond to the same state. So a pure state of the oscillator corresponds to a one-dimensional subspace of $L^2(\mathbb{R})$, the subspace generated by any one of the norm-one functions corresponding to it. Thus we can identify the pure states of the oscillator with one-dimensional projection operators on $L^2(\mathbb{R})$. The (general) states are so-called density operators, which are positive trace-class operators D of trace 1. Each such can be written $D = \sum w_i R_i$ where $w_1$, $w_2$, ... is a finite or countably infinite sequence of positive real numbers of sum 1, and the $R_i$ are pairwise orthogonal one-dimensional projections [B-C, §2.1(ii)]. For example, consider the first two eigenstates $\Psi_0$ and $\Psi_1$ of the oscillator. The linear combination $\Psi = (3/5)\Psi_0 + (4/5)\Psi_1$ is also a pure state of the oscillator because $\| \Psi \|^2 = (3/5)^2 + (4/5)^2 = 1$ but is not an eigenstate. If $R_0$ and $R_1$ are the projections on the subspaces generated by $\Psi_0$ and $\Psi_1$ respectively, then $D = (3/5)^2 R_0 + (4/5)^2 R_1$ is a density operator (a state), but not a projection (because $D^2 \neq D$), thus corresponds to no pure state.

As we have mentioned, the observables of the linear oscillator correspond to self-adjoint operators on $L^2(\mathbb{R})$; the special observables called questions, whose set is $\mathcal{L}$, correspond to the projection operators. Given a (general) state of our oscillator, which we may take as a density operator $D$, what then is the corresponding probability measure $m$? It is $m(Q) = $ trace $(QD) = $ trace $(DQ)$, $Q \in \mathcal{L}$. If $D$ corresponds to a pure state, then $D = P$, where $P$ is the one-dimensional projection on the subspace $\mathbb{C}\Psi$ generated by a norm-one function $\Psi \in L^2(\mathbb{R})$. In this case the measure $m$ corresponding to the pure state $D$ is $m(Q) = <Q(\Psi), \Psi>$, $Q \in \mathcal{L}$. The support of this pure state is the one-dimensional projection $P$ itself, $P \in \mathcal{L}$. Now observe the following two facts: (1) given any nonzero question $Q \in \mathcal{L}$, there is a pure state $m$ with $m(Q) = 1$; and (2) if $P$ is the support of a pure state $m$ and $n$ is any other state, with $n(P) = 1$, then $n = m$. As for (1), simply choose any norm-one function $\Psi$ in the subspace on which $Q$ projects. Then $<Q(\Psi), \Psi> = <\Psi, \Psi> = 1$. As for (2) it forces trace $(DP) = 1$, from which it follows readily that $D = P$. We shall take these two facts as axioms for the pure states.

**Axiom B.** (1) Given any nonzero question $a \in \mathcal{L}$, there is a pure state $m \in \mathcal{S}$ with $m(a) = 1$. (2) If $m$ is a pure state with support $a \in \mathcal{L}$, then $m$ is the only state, pure or not, with $m(a) = 1$.

Axiom B, in exactly this form, was proposed by MacLaren almost thirty years ago (1965) [Mac, Axiom $10'$ on p. 15]. Both Zieler and MacLaren are students of Mackey. Axiom B appears in [B-C, p. 165] as well.

From Axiom A and B together, we deduce

**Lemma 5.4.**

1) *Every nonzero $a \in \mathcal{L}$ is the support of some state $m \in \mathcal{S}$.*
2) *The support of a pure state is an atom in $\mathcal{L}$, and every atom is the support of a unique state which is necessarily pure.*
3) *$\mathcal{L}$ is a complete, atomistic orthomodular lattice.*

We shall outline the proof with specific reference to [B-C] for most of the details.

1. Given $0 \neq a \in \mathcal{L}$, Axiom B1 tells us that there is a state $n \in \mathcal{S}$ with $n(a) = 1$; let $W(a) = \{x \in \mathcal{L} : x = $ support$(x)$ for some $n \in \mathcal{S}$ with $n(a) = 1\}$. Manufacture, via Zorn's lemma, a maximal orthogonal set $\{x_i\}$ in $W(a)$. If $x_i = $ support $(n_i)$, then $x = \bigvee x_i$ is the support of $m = \sum t_i n_i$ where the $t_i$ are positive real numbers of sum 1. Then $x = a$, lest we contradict the maximality of $\{x_i\}$ [B-C, p. 298].

2. Suppose that $m$ is pure and $a \in \mathcal{L}$ its support. If $a$ were not an atom, we could write $a = b \vee c$ for orthogonal nonzero $b$ and $c$. By (1), b = support $(m_1)$ and c = support $(m_2)$ for some $m_1$, $m_2 \in \mathcal{S}$. If $n = (1/2)m_1 + (1/2)m_2$, then $n(a) = 1$. So, by B2, $n = m$, a contradiction [B-C, p. 301]. Conversely, suppose $a \in \mathcal{L}$ is an atom. By B1, there is a pure state $m$ with $m(a) = 1$. Let $b = \text{support}(m)$. Then $a \geq b$ because $m(a) = 1$. As $a$ is an atom, we must have $b = a$, so $a = \text{support}(m)$. If $n \in \mathcal{S}$ satisfies $n(a) = 1$, then $n = m$ by B2.

3. This readily established using (1) and (2) [B-C, §11. 6].

Next we turn to the superposition principle. If there is one feature characteristic of this strange theory, it is this principle.

The linear oscillator is an excellent example to use to explain the superposition principle. Recall our earlier description of the oscillator. The observable "total energy" E corresponds to the selfadjoint operator $H$ on $L^2(\mathbb{R})$, and the normalized eigenfunctions $\Psi_n$, $n = 0, 1, 2, ...$ of $H$ are the (pure) eigenstates of the oscillator. The eigenvalues $E_n$ of $H$, $H(\Psi_n) = E_n \Psi_n$, $n = 0, 1, 2, ...$, are the possible measured values of the oscillator's total energy. Consider the first two: $E_0 = (1/2)h\nu_0$, the zero point energy, and $E_1 = 3E_0$. If the oscillator is in state $\Psi_0$, then a measurement of its total energy will give the value $E_0$ with certainty. If in state $\Psi_1$, then $E_1$ every time. Now the principle of superposition of states asserts that $\Psi = (3/5)\Psi_0 + (4/5)\Psi_1$ is a possible state for the single linear oscillator. Suppose the oscillator is in that pure state. Then a measurement of its total energy will yield either the value $E_0$ or the value $E_1$, never any in-between value. You will get $E_0$ with probability $(3/5)^2$ and $E_1$ with probability $(4/5)^2$. Consider the case where you have made the measurement and have obtained the value $E_1$. Then must not the oscillator have been in the state $\Psi_1$ prior to that measurement? Quantum mechanics says no. Prior to the measurement, the oscillator was in the superposed state $\Psi = (3/5)\Psi_0 + (4/5)\Psi_1$, partly in the eigenstate $\Psi_0$, and partly in the eigenstate $\Psi_1$. But that is like saying that photon is partly polarized along the x-axis and partly polarized along the y-axis [D, §2 of Chapter 1], or that the electron is partly here and partly there [F, Chapter 6]. It's crazy, but it's quantum mechanics.

Using Lemma 5.4, we shall state the superposition principle in terms of the atoms of $\mathcal{L}$ [B-C, §§11.5 and 14.8].

**Definition 5.5.** The atom $c \in \mathcal{L}$ is a **superposition** of the atoms $a$ and $b$ if $c \neq a$, $c \neq b$, and $c \leq a \vee b$.

**Axiom C** (Principle of superposition of pure states).

1. Given two different pure state (atoms) $a$ and $b$, there is at least one other pure state $c$, $c \neq a$, and $c \neq b$ that is a superposition of $a$ and $b$.
2. If the pure state $c$ is a superposition of the distinct pure states $a$ and $b$, then $a$ is a superposition of $b$ and $c$.

Axiom C(1) implies that $\mathcal{L}$ is irreducible [B-C, Theorem 14.8.4], and Axiom C(2) implies that $\mathcal{L}$ has the covering property [B-C, Theorem 14.8.10]. (For definitions see §4.) Remarkably, quantum-theoretical properties involving pure states and their superposition correspond to familiar lattice-theoretic properties.

Putting this all together, Axioms A, B, and C imply that $\mathcal{L}$ is an irreducible, complete, orthomodular AC-lattice. The arguments of §§3 and 4 then show that $\mathcal{L}$ can be realized as the lattice of all closed subspaces of an orthomodular space $\{E, \mathcal{K}, <\cdot,\cdot>\}$ where $\mathcal{K}$ is a $*$-field, $E$ a left vector space over $\mathcal{K}$, and $<\cdot,\cdot>$ an orthomodular form on $E$; this fact was known [B-C, Chapter 21].

Our final axiom brings in the notion of the symmetry group of our physical system [J, §9.4]. The system is symbolized by the orthomodular lattice $\mathcal{L}$ of all closed subspaces of the orthomodular space $\{E, \mathcal{K}, <\cdot,\cdot>\}$. Let us deal with the associated projective geometry-with-polarity $\{\mathbb{P}(E), \perp\}$; $\mathcal{L}$ is the lattice of all closed subspaces of $\{\mathbb{P}(E), \perp\}$ (see §4). The symmetry group referred to is the group of all automorphisms of $\{\mathbb{P}(E), \perp\}$, which are the projective-geometry automorphisms $\varphi$ of $\mathbb{P}(E)$ that also satisfy $\varphi(a^\perp) = \varphi(a)^\perp$ for all points $a$ in $\mathbb{P}(E)$ [J, §9.4]. We focus attention on the subgroup of this group consisting of all those automorphisms represented by unitary operators in $E$. A **unitary operator** on $E$ is a bijective linear map $U$ of $E$ onto itself that preserves the form:

$$U(\alpha x + \beta y) = \alpha U(x) + \beta U(y) \quad \forall \alpha, \beta \in \mathcal{K}, \quad \forall x, y \in E$$
$$<U(x), U(y)> = <x, y> \quad \forall x, y \in E.$$

We frame our final hypothesis in terms of unitary operators.

**Axiom D** (ample unitary group). Given any two orthogonal pure states $a$, $b$ in $\mathcal{L}$, there is a unitary operator $U$ such that $U(a) = b$.

The unitary operator $U$ is a linear mapping of the vector space $E$ onto itself. It induces a projective-geometry automorphism of $\{\mathbb{P}(E), \perp\}$ which we are denoting by the same symbol $U$. When we write $U(a) = b$, we mean it in the latter sense, because the pure states $a$ and $b$ are points in $\{\mathbb{P}(E), \perp\}$. Now $a = \mathcal{K}e$ and $b = \mathcal{K}f$ for orthogonal nonzero vectors $e$ and $f$ in $E$, so for the operator $U$ we shall have $U(e) = \alpha f$, for some $\alpha \in \mathcal{K}$. Then $<e, e> = <U(e), U(e)> = <\alpha f, \alpha f> = \alpha <f, f> \alpha^*$. From this it follows that there exists in $E$ an infinite orthogonal sequence $\{e_i : i = 1, 2, ...\}$ such that $<e_i, e_i> = <e_j, e_j>$ for all $i, j$. Solèr's Theorem then applies, and we have

**Theorem 5.6.** *Let $\{\mathcal{L}, \mathcal{S}\}$ represent a system of questions and states satisfying the requirements of Definitions* 5.1 *and* 5.2*. If the system $\{\mathcal{L}, \mathcal{S}\}$ satisifes these four additional conditions*:
*Axiom* A = *regularity properties of the states (probability measures)*,
*Axiom* B = *existence and uniqueness of pure states*,
*Axiom* C = *principle of superposition of pure states*,
*Axiom* D = *ample unitary group*,
*then Axiom* VII *follows as a necessary consequence: The quantum logic $\mathcal{L}$ is isomorphic as an orthocomplemented partially ordered set to the orthocomplemented* *partially ordered set of all closed subspaces of a separable real, complex, or quaternionic Hilbert space.*

Department of Mathematics and Statistics, University of Massachusetts, Box 34515, Amherst, Massachusetts 01003-4515

*E-mail address*: `holland@math.umass.edu`